\documentclass[10pt]{article}

\usepackage[hidelinks=true]{hyperref}
\hypersetup{                    
    colorlinks=true,                
    breaklinks=true,                
    urlcolor= red,                 
    linkcolor= blue,                
    citecolor= blue                
    }

\usepackage[babel=true]{csquotes}
\usepackage{epigraph}
\usepackage{graphicx}
\usepackage{eurosym}
\usepackage{subcaption}
\usepackage{listings}
\usepackage{dsfont}
\usepackage{fancybox}
\usepackage[x11names,table,dvipsnames]{xcolor}
\usepackage{hyperref}
\usepackage[all]{hypcap}
\usepackage{wrapfig}
\usepackage{enumitem}
\usepackage{amsmath}
\usepackage{mathtools}
\usepackage[frak=mma]{mathalfa}
\usepackage{yfonts}
\usepackage{amsthm}
\usepackage{amsmath}
\usepackage{amssymb}
\usepackage{textcomp}
\usepackage{mathrsfs}
\usepackage[utf8]{inputenc}  
\usepackage[ruled,vlined]{algorithm2e}
\usepackage{indentfirst}
\usepackage[left=2cm,right=2cm,top=1cm,bottom=2cm]{geometry}
\usepackage[affil-it]{authblk}

\usepackage{pgfplotstable}
\usepackage{pgf,pgfarrows,pgfnodes,pgfautomata,pgfheaps,pgfshade}
\usepackage{colortbl}
\usepackage[english]{babel}
\usepackage{times}
\usepackage{graphicx}
\usepackage{overpic}
\usepackage{verbatim}
\usepackage{algorithmic}
\usepackage{multicol}

\usepackage{csvsimple}
\usepackage{multirow}
\usepackage{siunitx}
\usepackage{mathdots}

\usepackage{multirow}
\usepackage{stmaryrd}
\usepackage{float}

\usepackage{tikz}
\usepackage{etex}
\usepackage{pgfplotstable}
\usetikzlibrary{plotmarks,backgrounds}
\usepgfplotslibrary{units}
\usepgfplotslibrary{dateplot}
\usepackage{multimedia}
\usepackage{fp}
\usepackage{bbding}
\usetikzlibrary{shadows,spy,fixedpointarithmetic,patterns,intersections,arrows,shapes,decorations.text,decorations.pathmorphing,backgrounds,fit,positioning,shapes.symbols,chains}
\usepackage{siunitx}
\usepackage{booktabs,colortbl}
\usepackage{csquotes}
\usepackage{listings}

\usepackage{etoolbox}

\usepackage{arydshln}

\pgfplotstableread{
Year    Far           Near       Here           There          It
1024    1.7786e+02    264.033    1.125198e+01   7.741145e+01   28
2048    6.2690e+01    97.292     5.675948e+00   2.038585e+01   23
4096    1.9644e+01    35.701     8.044349e+00   9.729089e+00   20
8192    6.3323e+00    22.0793    1.730299e+01   6.047999e+00   27
}\datatableA
\pgfplotstableread{
Year    Far           Near       Here           There          It
1024    1.7786e+02    131.762    1.014010e+01   3.429036e+01   28
2048    6.2690e+01    53.828     6.202990e+00   1.752040e+01   28
4096    1.9644e+01    27.0714    5.095264e+00   8.638499e+00   23
8192    6.3323e+00    20.7763    6.963561e+00   4.791503e+00   23
}\datatableB
\pgfplotstableread{
Year Far           Near           It
1    1.291981e+01  2.345384e+00   17
2    5.970098e+00  0              0
3    7.532564e+01  8.595954e+00   16
4    1.637180e+01  5.671643e+00   25
}\datatableC
\pgfplotstableread{
Year Far           Near           It
1    8.151199e+01  2.261496e+01   29
2    1.460417e+01  0              0
3    7.780524e+01  0              0
4    1.745146e+01  0              0
}\datatableD
\pgfplotstableread{
Year    Far           Near       Here           There          It
256     6.4238e+01    117.73     2.778861e+00   1.581368e+01   13
512     6.3968e+01    112.16     3.346933e+00   1.992970e+01   18
1024    6.3215e+01    118.57     4.423136e+00   1.618159e+01   14
2048    5.9428e+01    117.58     6.908502e+00   2.133982e+01   17
4096    5.8143e+01    110.68     1.075236e+01   2.789055e+01   20
8192    5.4959e+01    116.63     1.991908e+01   2.364040e+01   17
}\datatableE
\pgfplotstableread{
Year    Far           Near       Here           There          It
256     2.9402e+01    111.345    9.392020e+00   2.570738e+01   29
512     2.9604e+01    111.518    9.959257e+00   2.798867e+01   28
1024    2.9434e+01    112.179    9.921512e+00   3.362715e+01   28
2048    2.9181e+01    112.232    1.004730e+01   3.373598e+01   28
4096    2.9798e+01    113.689    1.087214e+01   3.102074e+01   26
8192    2.9826e+01    113.807    1.326606e+01   3.066839e+01   25
}\datatableF

\makeatletter
\patchcmd{\epigraph}{\@epitext{#1}}{\itshape\@epitext{#1}}{}{}
\makeatother

\title{A Unified Framework for Double Sweep Methods for the Helmholtz Equation}
\author{{Nacime} Bouziani\textsuperscript{}\thanks{\noindent\textsuperscript{}Department of Mathematics, Imperial College London, London, SW7 2AZ, UK (n.bouziani18@imperial.ac.uk)
}\,,
{Fr\'ed\'eric} Nataf and Pierre-Henri Tournier \textsuperscript{}{\footnote{\textsuperscript{}
Laboratoire J.L. Lions, Universit\'e Pierre et Marie Curie, 4 place Jussieu, 75005 Paris, France, and ALPINES INRIA, Paris, France (frederic.nataf@sorbonne-universite.fr, pierre-henri.tournier@sorbonne-universite.fr)}} }
\date{}

\newcommand{\reels}{\mathbb{R}}

\newenvironment{remarque}{\textbf{\newline Remark : }}{}
\newtheorem{theorem}{Theorem}[section]    
\newtheorem{remark}{Remark}[section]    
\newtheorem{proposition}[theorem]{Proposition}
\newtheorem{madefinition}[theorem]{Definition} 

\newtheorem{lemma}[theorem]{Lemma}

\newcommand{\reviewerone}[1]{ { #1}}
\newcommand{\reviewertwo}[1]{ {  #1}}
\newcommand{\addedinrevision}[1]{ {  #1}}

\begin{document}


\maketitle
\setcounter{tocdepth}{5}
\tableofcontents

\begin{abstract}
	We consider sweeping domain decomposition preconditioners to solve the Helmholtz equation in the case of stripwise domain decomposition with or without overlaps. We unify their derivation and convergence studies by expressing them as Jacobi, Gauss-Seidel, and Symmetric Gauss-Seidel methods for different numbering of the unknowns. The proposed framework enables theoretical comparisons between the double sweep methods in~\cite{Nataf:1997:CRO,vion:2018:improved} and \reviewertwo{those} in~\cite{Stolk:2013:rapidly,vion:2014:double,Stolk:2017:improved}. Additionally, it facilitates the introduction of a new sweeping algorithm. We provide numerical test cases to assess the validity of the theoretical studies.
\end{abstract}

\section{Introduction}

Solving the Helmholtz equation numerically is a difficult task, especially when dealing with high-frequency regimes, heterogeneous media or reflecting boundary conditions. Over the last decades a lot of effort and progress has been made in developing efficient algorithms to solve the ill-conditioned linear system resulting from the Helmholtz operator's discretization. Domain decomposition methods (DDM) try to overcome these difficulties. They are hybrid methods that combine direct solvers in subdomains and iterative matching of the solutions across the subdomains. The original domain decomposition method introduced by Schwarz~\cite{Schwarz:1870:UGA} only works for overlapping domain decomposition. P. L. Lions~\cite{Lions:1988:SAM} introduced a new variant of this algorithm where the Dirichlet interface conditions are replaced by Robin interface conditions, his method can be applied to both overlapping and nonoverlapping subdomains. He showed convergence for the elliptic case for a non overlapping domain decomposition. The proof was extended by Despr\'es~\cite{Despres:1991:DDM} to the Helmholtz equation and later on to the time-harmonic Maxwell equations~\cite{Despres:1992:ADD}. More recently, sweeping-type domain decomposition methods  have been made popular due to their capability to achieve nearly-linear asymptotic complexity. A sweeping algorithm was first proposed and analyzed  in~\cite{Nataf:1997:CRO} for convection-diffusion operators. Sweeping approaches for Helmholtz problems have recently seen their interest renewed as a preconditioner to speed up the convergence of the solver: the double sweep preconditioners of Vion and Geuzaine for non overlapping decomposition with high order interface conditions  \cite{vion:2014:double,vion:2018:improved}, the PML-based sweep method of Stolk \cite{Stolk:2013:rapidly,Stolk:2017:improved} for overlapping decomposition and the polarized traces method of Zepeda-N\'u\~nez and Demanet \cite{zepeda:2016:method}. \addedinrevision{Shifted Laplace preconditioners~\cite{Magolu:2001:incomplete,Erlangga:2004:Class} for wave propagation problems may be used in conjonction with DD methods, see~\cite{Graham:2017:ddHelmholtz,Bonazzoli:2019:domain} and references therein.} There also exist sweeping-type methods that are not domain decomposition based methods, such as the sweeping PML preconditioner of Engquist and Ying \cite{Engquist:2011:sweeping,Engquist:2011:SweepingHierarchical}, the source transfer method~\cite{Chen:2013:source} , see~\cite{gander:2019:class} for a complete panorama and relations between these methods. \addedinrevision{Another fruitful line of research consists in using time domain solvers. In~\cite{Mur:1993:finite}, the problem is solved on a time interval long enough so that the solution is periodic in time. A more efficient approach consists in searching for a time periodic solution of the time dependent wave equation using controllability techniques, see  \cite{Bristeau:1998:controllability} for the seminal paper and for recent works, see e.g. \cite{Tang:2022:Fully} and references herein.}\\

Our main contributions are as follows:
\begin{itemize}
	\item New formulation of \cite{Stolk:2013:rapidly,vion:2014:double} which results in \reviewertwo{the new variant} introduced in \S~\ref{sec:algosLR}
	\item Unified convergence analysis for the three above mentioned algorithms which enables a comparison with the algorithms proposed in \cite{Despres:1992:ADD, Stolk:2013:rapidly,Nataf:1997:CRO,vion:2018:improved}, see Table~\ref{tab:algocvrate}.
	\item \addedinrevision{Numerical comparisons between the different sweeping algorithms that align with the theoretical results.} 
\end{itemize}

We first state the problem in \S~\ref{sec:problem_setting}. Then we explain in~\S~\ref{sec:substructuring} how to substructure the problem in terms of interface unknowns and how to apply classical linear algebra preconditioners (Jacobi,Gauss-Seidel and Symmetric Gauss-Seidel) to two different unknown numberings. After these preparatory tools have been introduced, we present and analyze in~\S~\ref{sec:convergence_rates} the convergence of six sweeping algorithms. Numerical results are shown in~\S~\ref{sec:resnum}.

\section{Statement of the problem and two classical algorithms} 
\label{sec:problem_setting}

\begin{figure}[htp]
  \centering
  \includegraphics[width=0.6\textwidth]{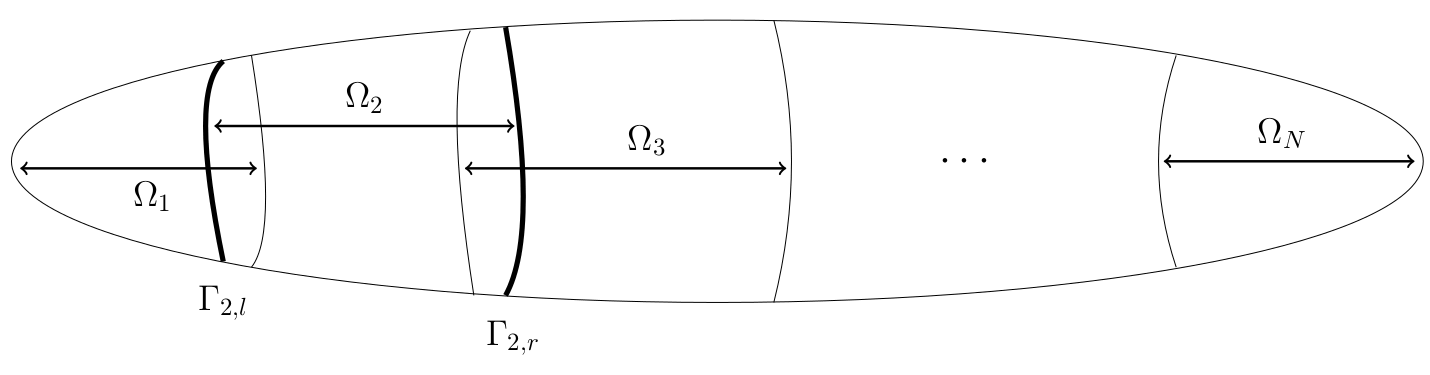}
  \caption{Decomposition into vertical strips}
  \label{fig:decomposition}
\end{figure}

\begin{figure}
\centering
\includegraphics[width=0.6\textwidth]{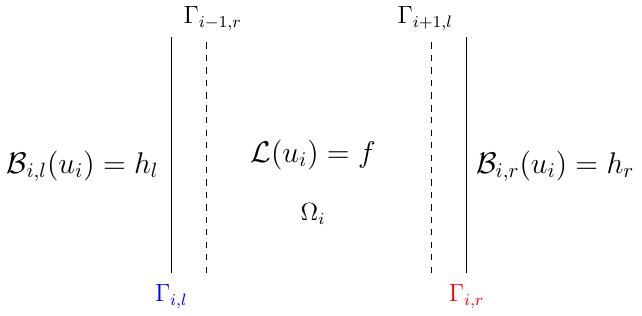}
\caption{Local problem on the subdomain $\Omega_{i}$
     }
\label{fig:subdomain}
\end{figure}

We consider the Helmholtz equation in a bounded domain $\Omega \subset \reels^{2}$ with frequency $\omega$, velocity $c$ and wavenumber $k$ defined by $k = \omega/c$:
\begin{equation}
\label{helmholtz_equation}
\begin{array}{l}
\left(- k^2 - \Delta\right) u = f \text{ in }\Omega\\
	+ \text{ appropriate boundary conditions on }\partial\Omega\,.
\end{array}
\end{equation}
We consider a layered decomposition of $\Omega$ into $N$ slices $(\Omega_i)_{1\le i \le N}$ with or without overlap, see Figure~\ref{fig:decomposition}. The boundary $\partial\Omega_i \setminus \partial\Omega$ is written as the disjoint union of $\Gamma_{i,l}$ and $\Gamma_{i,r}$ where $\Gamma_{i,l}$ is on the left of $\Omega_i$ and $\Gamma_{i,r}$ is on its right ($\reviewerone{\Gamma}_{1,l}=\emptyset$ and $\reviewerone{\Gamma}_{N,r}=\emptyset$) (see Fig.~\ref{fig:subdomain}). The outward normal from $\Omega_i$ on $\Gamma_{i,l}$ (resp. $\Gamma_{i,r}$) is denoted by $\vec n_{i,l}$ (resp. $\vec n_{i,r}$). 
The problem~\eqref{helmholtz_equation} can be solved iteratively using a domain decomposition method where we solve locally on each subdomain $\Omega_i$ the equation \eqref{helmholtz_equation} with appropriate boundary conditions on the physical boundaries and interfaces~\cite{Despres:1991:DDM}. The method writes:\\
 Solve in parallel:
\begin{equation}
\label{jacobi_algo_helmholtz}
\left\{
\begin{aligned}
\left(- k^2 - \Delta\right) u_{i}^{n+1} &= f\ \text{ in } \Omega_{i},\ 1\le i\le N \\
\mathcal{B}_{i,l}\left(u_{i}^{n+1}\right) &= \mathcal{B}_{i,l}\left(u_{i-1}^{n}\right) \ \text{on} \ \Gamma_{i,l},\ 2\le i\le N \\
\mathcal{B}_{i,r}\left(u_{i}^{n+1}\right) &= \mathcal{B}_{i,r}\left(u_{i+1}^{n}\right) \ \text{on} \ \Gamma_{i,r},\ 1\le i\le N-1 \\
&+ \text{ appropriate boundary conditions on }\partial\Omega\cap \partial\Omega_i\,,
\end{aligned}
\right.
\end{equation}
where $\mathcal{B}_{i,l}$ and $ \mathcal{B}_{i,r}$ are the interface conditions. Here $I$ denotes the square root of $-1$ ($I^2=-1$). For sake of simplicity, we consider here either zeroth-order ABC (ABC0):
\begin{equation}
	\label{eq:IBC0} 
\left\{
\begin{aligned}
\mathcal{B}_{i,l} &= \partial_{\vec{n}_{i,l}} + Ik  \\
\mathcal{B}_{i,r} &= \partial_{\vec{n}_{i,r}} + Ik \,,
\end{aligned}
\right.
\end{equation}
exact ABC as interface conditions or subdomains coupling via Perfectly Matched Layer (PML) as in~\cite{Stolk:2013:rapidly}. In practice, when used as truncation conditions on artificial boundaries, ABC0 yields high non physical reflection of the order of 10\% of the incoming wave whereas by definition exact ABC leads to no reflection at all. ABC0 boundary conditions are easy and cheap to use whereas exact ABC, always defined in theory, are sometimes impossible to use in practice (e.g. variable coefficients problems). As a result, there is room for compromise and a great deal of literature has been devoted to introduce various high order interface conditions. The most notable techniques are based on partial differential operators (see~\cite{Enquist:1977:ABC} and \cite{Bayliss:80:RBC}) or on PML (see~\cite{PML} and \cite{Chew:1994:3PM}), see~\cite{nataf:2013:absorbing} as well for a gentle introduction to this question. When used as truncation conditions, the final accuracy of the computation does depend on the choice of the ABC. 

But let us stress that here, ABCs are used as interface conditions in domain decomposition methods so that the final accuracy of the computed result is not impacted by the choice of the ABC, only the iteration counts to solution are impacted. High-order ABC tailored to domain decomposition improve the iteration counts with respect to ABC0, see e.g. ~\cite{Gander:2001:OSH,Antoine:12:QON}. A remarkable super convergence result noticed in~\cite{Hagstrom:1988:NED} for the two subdomain case and in~\cite{Nataf:1994:OIC} for a decomposition into $N$ strips (see Fig.~\ref{fig:decomposition}) is that the use of exact ABCs as interface conditions yields convergence in a number of iterations equal to the number of subdomains. Since the solution in a subdomain depends on the value of the right hand side everywhere and that in algorithm~\eqref{jacobi_algo_helmholtz} a subdomain receives data only from its neighbors it is not possible to achieve convergence in less than $N$ iterations. When looking at the proof in~\cite{Nataf:1994:OIC}, it appears that somehow the correct information flows from the extreme subdomain labelled $1$ to the right and at the same time from subdomain labelled $N$ to the left.

This motivated the search for algorithms which would sweep over the subdomains to reach convergence in one iteration consisting of a double sweep. In the sequel we will consider the double sweep algorithms introduced in~\cite{Nataf:1997:CRO,vion:2018:improved,Stolk:2013:rapidly,vion:2014:double} since they converge in one double sweep if exact ABCs are used as interface conditions. Note that these algorithms were named double sweep algorithms which could be confusing. But as we shall see they are actually not the same and have different convergence rates when implemented with non exact ABCs. We start with the double sweep algorithm introduced in~\cite{Nataf:1997:CRO}. It consists in double sweeps over the subdomains:\\
{\bf left to right sweep:}
\begin{equation}
\label{leftsweep_algo_helmholtz}
\left\{
\begin{aligned}
\left(- k^2 - \Delta\right) u_{i}^{n+1/2} &= f\ \text{ in } \Omega_{i},\ 1\le i\le N \\
\mathcal{B}_{i,l}\left(u_{i}^{n+1/2}\right) &= \mathcal{B}_{i,l}\left(u_{i-1}^{n+1/2}\right) \ \text{on} \ \Gamma_{i,l},\ 2\le i\le N \\
\mathcal{B}_{i,r}\left(u_{i}^{n+1/2}\right) &= \mathcal{B}_{i,r}\left(u_{i+1}^{n}\right) \ \text{on} \ \Gamma_{i,r},\ 1\le i\le N-1 \\
&+ \text{ appropriate boundary conditions on }\partial\Omega\cap \partial\Omega_i\,,
\end{aligned}
\right.
\end{equation}
{\bf then, right to left sweep:}
\begin{equation}
\label{rightsweep_algo_helmholtz}
\left\{
\begin{aligned}
\left(- k^2 - \Delta\right) u_{i}^{n+1} &= f\ \text{ in } \Omega_{i},\ 1\le i\le N \\
\mathcal{B}_{i,l}\left(u_{i}^{n+1}\right) &= \mathcal{B}_{i,l}\left(u_{i-1}^{n+1/2}\right) \ \text{on} \ \Gamma_{i,l},\ 2\le i\le N \\
\mathcal{B}_{i,r}\left(u_{i}^{n+1}\right) &= \mathcal{B}_{i,r}\left(u_{i+1}^{n+1}\right) \ \text{on} \ \Gamma_{i,r},\ 1\le i\le N-1 \\
&+ \text{ appropriate boundary conditions on }\partial\Omega\cap \partial\Omega_i\,.
\end{aligned}
\right.
\end{equation}
In \reviewerone{\S~\ref{sec:algosSW}, we show} it can be seen as a Symmetric Gauss-Seidel (SGS) version of the Jacobi algorithm~\eqref{jacobi_algo_helmholtz}. 
Its study is made easier in its substructured formulation which is moreover needed to introduce the double sweep algorithms of~\cite{Stolk:2013:rapidly,vion:2014:double}. We devote the next section to substructuring. 


\section{Substructuring}
	\label{sec:substructuring}

In this section, we introduce the substructured problem related to algorithm~\eqref{jacobi_algo_helmholtz}. This will be the basis for the unified framework of DD sweeping methods.  Substructuring consists in reformulating the iterative method considering only surface unknowns on the interfaces:
\begin{equation}
\label{hl_hr_reformulation}
\left\{
\begin{aligned}
h^{n}_{i,l} &\vcentcolon = \mathcal{B}_{i,l}\left(u^{n}_{i}\right), \ \text{ on }\ \Gamma_{i,l} \ \ \text{ for } \ 2\le i\le N \\
h^{n}_{i,r} &\vcentcolon =  \mathcal{B}_{i,r}\left(u^{n}_{i}\right), \ \text{ on }\ \Gamma_{i,r} \ \ \text{ for } \ 1\le i\le N-1\,.
\end{aligned}
\right.
\end{equation}
Considering the global vector $h^n$ containing the local unknowns $(h_{i,l}^{n})_{2\le i \le N}$ and $(h_{i,r}^{n})_{1\le i \le N-1}$, we can reformulate the additive Schwarz method~\eqref{jacobi_algo_helmholtz} as a Jacobi algorithm on $h^n$:
\begin{equation}
\label{global_jacobi}
h^{n+1} \vcentcolon= \mathscr{T}(h^{n}) + G\,,
\end{equation}
where the iteration operator $\mathscr{T}$ can be written in the form of an operator valued matrix and $G$ refers to the contribution of the right-hand side $f$, see~\cite{Nataf:1994:OIC}. The above equation is what is called a substructured formulation of the volumic algorithm~\eqref{jacobi_algo_helmholtz}. Taking the limit as $n$ tends to infinity, we see that we look for a vector $h$ such that,
\begin{equation}
\label{global_substructuring}
\left( Id - \mathscr{T}\right)(h) = G\,.
\end{equation}
Equation~\eqref{global_substructuring} is what is called the substructured formulation of the domain decomposition problem. In order to define more precisely the operator $\mathscr{T}$,
 we introduce for each subdomain an operator $S_{i}$ which takes three arguments, two surface functions $h_{l}$ and $ h_{r}$ and a volume function $f$ and maps them to the local solution $v$:
 \begin{equation}
	 \label{eq:Si}
 S_i(h_{l}, h_{r}, f) \vcentcolon= v 	
 \end{equation}
 where $v:\,\Omega_i \mapsto \mathbb{C}$ satisfies:
\begin{equation}
\label{operator_S_helmholtz}
\left\{
\begin{aligned}
\left(- k^2 - \Delta\right) v &= f \ \ \text{in} \ \Omega_{i} \\
\mathcal{B}_{i,l}(v) &= h_{l} \ \text{on} \ \Gamma_{i,l} \quad \left(2\le i \le N\right)\\
\mathcal{B}_{i,r}(v) &= h_{r} \ \text{on} \ \Gamma_{i,r} \quad \left(1\le i \le N-1\right)\\
&+ \text{ appropriate boundary conditions on }\partial\Omega\cap \partial\Omega_i\,,
\end{aligned}
\right.
\end{equation}
for  $1 < i < N$. For $i=1$, the definition of $S_1$ is similar except that it takes only the two arguments $(h_r,f)$ since domain $\Omega_1$ has no left interface and similarly operator $S_N$ takes only the two arguments $(h_l,f)$ since domain $\Omega_N$ has no right interface. {\bf As of now, for sake of simplicity and by abuse of notation, $S_1(h_l,\ h_r,\ f )$ $(\text{resp. } S_N(h_l,\ h_r,\ f))$  will refer to $S_1(h_r,f)$ $(\text{resp. } S_N(h_l,\ f))$. }

Next, we introduce the surface right hand-side $G(f)$ by
\begin{equation}
	\label{eq:Gglobal}
	\begin{array}{rcl}
G_{i,l}(f) &\vcentcolon =& \mathcal{B}_{i, l}(S_{i-1}(0 , 0 , \ f)),\ \ 2\le i \le N \\
G_{i,r}(f)  &\vcentcolon =&  \mathcal{B}_{i, r}(S_{i+1}(0 , 0 , \ f)),\ \ 1\le i \le N-1\,.
	\end{array}
\end{equation}
and the substructured operator $\mathscr{T}$ by:
\begin{equation}
	\begin{array}{rcl}
\mathscr{T}(h)_{i+1,l} &\vcentcolon =& \mathcal{B}_{i+1, l}(S_i(h_{i,l},\ h_{i,r}, 0)),\ \ 1\le i \le N-1 \\
\mathscr{T}(h)_{i-1,r}  &\vcentcolon =&  \mathcal{B}_{i-1, r}(S_i(h_{i,l},\ h_{i,r}, 0)),\ \ 2\le i \le N\,.
	\end{array}
\end{equation}
The operator $\mathscr{T}$ has thus the following possibly non zero entries:
\begin{equation}
	\label{eq:Tentries}
	\begin{array}{rcl}
\mathscr{T}_{(i+1,l)(i,l)} &\vcentcolon =& \mathcal{B}_{i+1, l}(S_i(\cdot,0, 0)),\ \ 1\le i \le N-1 \\
\mathscr{T}_{(i+1,l)(i,r)} &\vcentcolon =& \mathcal{B}_{i+1, l}(S_i(0, \cdot , 0)),\ \ 1\le i \le N-1 \\
\mathscr{T}_{(i-1,r)(i,r)}  &\vcentcolon =&  \mathcal{B}_{i-1, r}(S_i(0,\cdot , 0)),\ \ 2\le i \le N \\
\mathscr{T}_{(i-1,r)(i,l)}  &\vcentcolon =&  \mathcal{B}_{i-1, r}(S_i(\cdot , 0, 0)),\ \ 2\le i \le N \,.
	\end{array}
\end{equation}

\subsection{Double sweep algorithm of \cite{vion:2014:double}}
	\label{sec:stolkviongeuzaine}
	
The rationale behind the algorithm in~\cite{vion:2014:double} is that for exact ABC used as interface conditions, we have a decoupling of left and right interface unknowns since  	
\begin{equation}
	\label{eq:TEABC}
	\mathscr{T}^{EABC}_{(i+1,l)(i,r)} \equiv 0 \text{ and } \mathscr{T}^{EABC}_{(i-1,r)(i,l)} \equiv 0 \,.
\end{equation}
Then, the operator $\mathscr{T}^{EABC}$ is nilpotent of order $N-1$. This is related to the convergence in $N$ iterations of algorithm~\eqref{jacobi_algo_helmholtz} with exact ABC as interface conditions. In practice, ABC and even PML truncation techniques are not perfect and the nilpotency effect is lost. In order to force it, a new operator is introduced:
\begin{equation}
	\label{eq:TsdsDefinition}
	\begin{array}{rcl}
{\mathscr{T}_{BJ}}(h)_{i+1,l} & :=& \mathcal{B}_{i+1, l}(S_i(h_{i,l} , 0 , 0)),\ \ 1\le i \le N-1 \\
{\mathscr{T}_{BJ}}(h)_{i-1,r}  & :=&  \mathcal{B}_{i-1, r}(S_i( 0 ,\ h_{i,r}, 0)),\ \ 2\le i \le N
	\end{array}
\end{equation}
which is by construction nilpotent of order $N-1$ even with non exact ABCs as interface conditions. A fixed point method based on this reads:
\begin{equation}
( Id - \mathscr{T}_{BJ})(h^{n+1}) =  \left(\mathscr{T} - \mathscr{T}_{BJ}\right)(h^n) + G.
\end{equation}
More efficiently, the operator $I - {\mathscr{T}_{BJ}}$ can then be used as a preconditioner in order to solve the substructured problem~\eqref{global_substructuring}. For instance the left preconditioned system 
\begin{equation}
\label{precondglobal_substructuring}
( Id - {\mathscr{T}_{BJ}}\, )^{-1} \left( Id - \mathscr{T}\right)(h) = ( Id - {\mathscr{T}_{BJ}}\, )^{-1}G
\end{equation}
can be solved by a Krylov type method. A closer look at the operator $I-{\mathscr{T}_{BJ}}$ shows that inverting it can be made by two concurrent sweeps, hence the name double sweep.\\ 

\medskip

In the sequel, algorithms other than \eqref{global_jacobi} or \eqref{precondglobal_substructuring} will be introduced as classical Jacobi, Gauss-Seidel, symmetric Gauss-Seidel applied to the substructured problem~\eqref{global_substructuring} with two different numberings of the interfaces. This will give a unified view to the methods considered in~\cite{Despres:1990:DDP,Nataf:1997:CRO,Stolk:2013:rapidly,vion:2014:double,vion:2018:improved} and it will enable the introduction of \reviewertwo{a} new algorithm in \S~\ref{sec:algosLR}.

\subsection{Three-subdomain case} 
\label{sub:threesubdomaincase}

For sake of clarity and to give a taste of the general case, we  start with a three-domain decomposition of the whole plane $\reels^2$, see Figure~\ref{fig:threesubdomain}.

\subsubsection{Subdomain wise (SW) ordering}
 The vector $h$ has four components which are interfaces functions living respectively on $\Gamma_{1,r}$, $\Gamma_{2,l}$, $\Gamma_{2,r}$ and $\Gamma_{3,l}$, see Figure~\ref{fig:threesubdomain}:
\[
	h := (h_{1,r}\ h_{2,l}\ h_{2,r}\ h_{3,l})^T\,.
\]
 This is the natural geometric ordering that we will refer to as the subdomain wise numbering (SW). Then for arbitrary interface conditions $\mathcal{B}_{i,l \text{ or } r}$, the sparsity pattern of the substructured problem is:
\begin{equation}
	\label{eq:T3sdsdordering}
	(I-\mathscr{T}_{SW})(h_{SW}) = \left(
\begin{array}{cccc}	
	I & X & X & 0 \\
	X & I & 0 & 0 \\ 
	0 & 0 & I & X \\
	0 & X & X & I 
\end{array}	
\right)\,
\left(\begin{array}{l}
  h_{1,r}\\
  h_{2,l}\\
  h_{2,r}\\
  h_{3,l}
\end{array}
\right)\,,
\end{equation}
where $X$ denotes a possibly non zero entry.

\subsubsection{Left-Right (LR) ordering}

Following~\cite{Nataf:1997:CRO}, we consider now an ordering where left interfaces are numbered first and then the right interfaces in the reverse order. In our three-subdomain case, we define;
\[
	h := (h_{2,l}\ h_{3,l}\ h_{2,r}\ h_{1,r})^T\,.
\]
For arbitrary interface conditions $\mathcal{B}_{i,l \text{ or } r}$, the sparsity pattern of the substructured operator is:
\begin{equation}
	\label{eq:T3sd}
	(I-\mathscr{T}_{LR})(h_{LR}) = \left(
\begin{array}{cc;{2pt/2pt}cc}	
	I & 0 & 0 & X \\
	X & I & X & 0 \\[.5em] \hdashline[2pt/2pt]
	0 & X & I & 0 \\
	X & 0 & X & I 
\end{array}	
\right)\,
\left(\begin{array}{l}
  h_{2,l}\\
  h_{3,l}\\[.5em] \hdashline[2pt/2pt]
  h_{2,r}\\
  h_{1,r}
\end{array}
\right)\,.\end{equation}

\subsubsection{Dependence of the preconditioners on ordering}
  \label{sec:threesditerative}

Of course, systems~\eqref{eq:T3sdsdordering} and~\eqref{eq:T3sd} are strictly equivalent. But as a consequence of the different numberings, the approximate inverses obtained from Jacobi or Gauss-Seidel type mechanisms may not be the same. As an example, consider the Gauss-Seidel preconditioner for the SW numbering:
\begin{equation}
	\label{eq:T3sdsdorderingBGS}
	(I-\mathscr{T}_{SW,GS})(h_{SW}) := \left(
\begin{array}{cccc}	
	I & 0 & 0 & 0 \\
	\mathscr{T}_{(2,l)(1,r)} & I & 0 & 0 \\ 
	0 & 0 & I & 0 \\
	0 & \mathscr{T}_{(3,l)(2,l)} & \mathscr{T}_{(3,l)(2,r)} & I 
\end{array}	
\right)\,
\left(\begin{array}{l}
  h_{1,r}\\
  h_{2,l}\\
  h_{2,r}\\
  h_{3,l}
\end{array}
\right)\,,
\end{equation}
and the Gauss-Seidel preconditioner for the LR numbering \reviewertwo{introduced in~\cite{Stolk:2013:rapidly,Stolk:2017:improved}}:
\begin{equation}
	\label{eq:T3sdBGS}
	(I-\mathscr{T}_{LR,GS})(h_{LR}) := \left(
\begin{array}{cc;{2pt/2pt}cc}	
	I & 0 & 0 & 0 \\
	\mathscr{T}_{(3,l)(2,l)} & I & 0 & 0 \\[.5em] \hdashline[2pt/2pt]
	0 & \mathscr{T}_{(2,r)(3,l)} & I & 0 \\
	\mathscr{T}_{(1,r)(2,l)} & 0 & \mathscr{T}_{(1,r)(2,r)} & I 
\end{array}	
\right)\,
\left(\begin{array}{l}
  h_{2,l}\\
  h_{3,l}\\[.5em] \hdashline[2pt/2pt]
  h_{2,r}\\
  h_{1,r}
\end{array}
\right)\,.\end{equation}
Even up to a reordering these two preconditioners are different since the entries of $\mathscr{T}$ which are kept in these two approximations are not the same. In \S~\ref{sec:convergence_rates}, we will develop this analysis for other preconditioners with also an arbitrary number of subdomains.\\

Another important fact already mentioned above, see eq.~\eqref{eq:TEABC}, is that if the interface conditions are exact absorbing boundary conditions (EABC), four additional entries cancel in~\eqref{eq:T3sdsdordering} or \eqref{eq:T3sd} so that $\mathscr{T}_{SW}$ and $\mathscr{T}_{LR}$ are actually nilpotent operators of order 2. Indeed the entry $\mathscr{T}_{(1,r)(2,l)}$ is the operator that maps $h_{2,l}$ to $\mathcal{B}_{1,r}(v_{2,l})$  where $v_{2,l}$ satisfies:
\begin{equation}
	\label{eq:AzeroEABC}
\begin{array}{l}
	\mathcal{L}(v_{2,l}) = 0 \text{ in }\Omega_2\,, \\
	\mathcal{B}_{2,r}(v_{2,l}) = 0 \text{ on }\Gamma_{2,r}\,, \\
	\mathcal{B}_{2,l}(v_{2,l}) = h_{2,l} \text{ on }\Gamma_{2,l}\,.
\end{array}
\end{equation}
Since $\mathcal{B}_{2,r}$ is an EABC, $v_{2,l}$ can be seen as the restriction of a harmonic function defined on the half plane at the right of $\Gamma_{2,l}$ so that $\mathcal{B}_{1,r}(v_{2,l})=0$ since $\mathcal{B}_{1,r}$ is an EABC as well. That is, $\mathscr{T}_{(1,r)(2,l)}=0$. In the same way, we can prove three other cancellations, namely: $\mathcal{T}_{(2,r)(3,l)}=0$ and $\mathcal{T}_{(i+1,l)(i,r)}=0$ for $i=1,2$. The only entries left are $\mathscr{T}_{(i,r)(i+1,r)}$ and $\mathscr{T}_{(i+2,l)(i+1,l)}$ for $i=1$. Let us denote by $\mathscr{T}^{EABC}_{SW}$ the operator $\mathscr{T}_{SW}$ when the interface conditions are EABC. We have just proved that its sparsity pattern is then:
\begin{equation}
	\label{eq:TEABC3sdsdordering}
	\mathscr{T}^{EABC}_{SW} = \left(
\begin{array}{cccc}	
	0 & 0 & X & 0 \\
	0 & 0 & 0 & 0 \\ 
	0 & 0 & 0 & 0 \\
	0 & X & 0 & 0 
\end{array}	
\right)\,.	
\end{equation}
As for the left-right numbering, let us denote by $\mathscr{T}^{EABC}_{LR}$ the operator $\mathscr{T}_{LR}$ when the interface conditions are EABC so that with four entries cancelled, its sparsity pattern is then:
\begin{equation}
	\label{eq:TEABC3sd}
	\mathscr{T}^{EABC}_{LR} = \left(
\begin{array}{cc;{2pt/2pt}cc}	
	0 & 0 & 0 & 0 \\
	X & 0 & 0 & 0 \\ \hdashline[2pt/2pt]
	0 & 0 & 0 & 0 \\
	0 & 0 & X & 0 
\end{array}	
\right)\,.
\end{equation}
Another way to see that numbering impacts the preconditioners is to remark that even for EABCs, the Gauss-Seidel preconditioner for the left-right numbering is exact:
\[
	(I-\mathscr{T}^{EABC}_{LR,GS})=(I-\mathscr{T}^{EABC}_{LR})\,,
\]
which is not the case for the  Gauss-Seidel preconditioner for the subdomain wise numbering
\[
	(I-\mathscr{T}^{EABC}_{SW,GS}) \neq (I-\mathscr{T}^{EABC}_{SW})\,. 
\]

\begin{figure}
\centering
\begin{tikzpicture}
	\def\deb{5}	
	\def\dx{1}
	\draw (\deb ,0) -- ++(0,4)       node[above] {$\Gamma_{2,l}$};
	\draw (\deb+\dx,0) -- ++(0,4)    node[above] {$\Gamma_{1,r}$};
	\draw (\deb+7,0) -- ++(0,4)      node[above] {$\Gamma_{3,l}$};
	\draw (\deb+7+\dx,0) -- ++(0,4)  node[above] {$\Gamma_{2,r}$};
	\draw[{stealth[scale=2]}-{stealth[scale=2]}] (\deb ,3) -- ++(7+\dx,0) 
	      node[above , midway] {$\Omega_{2}$};
	\draw[<-] (\deb+\dx ,2) -- ++(-5,0)       node[above , midway] {$\Omega_{1}$};
	\draw[<-] (\deb+7 ,2) -- ++(6,0)       node[above , midway] {$\Omega_{3}$};
\end{tikzpicture}
\caption{Decomposition of the plane into three subdomains}
\label{fig:threesubdomain}
\end{figure}


 \section{Analysis for the N-subdomain case} 
 \label{sec:convergence_rates}
We come back to the case of a stripwise decomposition into $N$ subdomains as in Figures~\ref{fig:decomposition} or~\ref{fig:subdomain}. As in the three-subdomain case, we consider two numberings for the substructured system~\eqref{global_substructuring}.

\subsection{Subdomain wise (SW) numbering}

The most natural numbering is to order the interface unknowns subdomain-wise (SW), i.e. one subdomain after the other. It yields the following substructured system:
\begin{equation}
	\label{eq:substructSW}
 (I-{\mathscr T}_{SW})(H_{SW})=
\left[ 
\begin{array}{cccccccc}
   I&-{\mathscr T}_{(1,r)(2,l)} &-{\mathscr T}_{(1,r)(2,r)} & 0 &\phantom{{-\mathscr T}}&\phantom{{-\mathscr T}} &\phantom{{-\mathscr T}} &\phantom{{-\mathscr T}} \\
   -{\mathscr T}_{(2,l)(1,r)}& I & 0 & 0 &\ddots &\ &\ &\  \\
   0 & 0 & I &-{\mathscr T}_{(2,r)(3,l)} &-{\mathscr T}_{(2,r)(3,r)} & 0 &\ &\  \\
   \ &{\mathscr T}_{(3,l)(2,l)} &{\mathscr T}_{(3,l)(2,r)}& I & 0 & 0 &\ &\  \\
   \ & \ & 0 & 0 &\ddots &\ddots &\ddots &\  \\
   \ &\ &\ &-{\mathscr T}_{(4,l)(3,l)}  &\ddots &\ddots &\ddots &\ddots  \\
   \ &\ &\ &\ &\ddots &\ddots &\ddots &\ddots  \\
   \ &\ &\ &\ &\ &\ddots\ &\ddots\ &I   
\end{array}
\right]\ 
\left[\begin{array}{c} h_{1,r}\\[.0em] h_{2,l} \\[.0em] h_{2,r} \\[.0em] h_{3,l}\\[.2em]
h_{3,r}\\[.0em] h_{4,l}\\[.4em] \vdots \\[.4em] h_{N,l}
\end{array}\right]\,.
\end{equation}

\subsection{Left-Right (LR) numbering}
The following left-right (LR) numbering, where left interface unknowns are numbered first and then the right interface unknowns, leads to the following substructured system:
\begin{equation}
	\label{eq:substructLR}
 (I-{\mathscr T}_{LR})(H_{LR})=
\left[ \begin{array}{c;{2pt/2pt}c}
\begin{array}{cccc}
   I& \ & \ & 0 \\
   -{\mathscr T}_{(3,l)(2,l)}&\ddots& \ &\ \\
   \ &\ddots&\ddots&\ \\
   0&\ &-{\mathscr T}_{(N,l)(N-1,l)}&I 
   \end{array}&
\begin{array}{cccc}
  0&\ &\ &-{\mathscr T}_{(2,l)(1,r)}\\
   \ &\ &\iddots &\ \\
   \ &\iddots &\ &\ \\
   -{\mathscr T}_{(N,l)(N-1,r)}&\ &\ &0 
   \end{array}\\ \hdashline[2pt/2pt]
\begin{array}{cccc}
  0&\ &\ &-{\mathscr T}_{(N-1,r)(N,l)}\\
   \ &\ &\iddots &\ \\
   \ &\iddots &\ &\ \\
   -{\mathscr T}_{(1,r)(2,l)}&\ &\ &0 
   \end{array} &  
   \begin{array}{cccc}
   I&\ \ &\ \ &0\\[.5em]
   -{\mathscr T}_{(N-2,r)(N-1,r)}&\ddots&\ &\ \\[.5em]
   \ &\ddots&\ddots&\ \ \\[.5em]
   0&\  &-{\mathscr T}_{(1,r)(2,r)}&I 
   \end{array}
\end{array}\right]\ 
\left[\begin{array}{c} h_{2,l}\\ \vdots \\ \vdots \\ h_{N,l}\\ \hdashline[2pt/2pt] \\
h_{N-1,r}\\ \vdots\\ \vdots \\ h_{1,r}
\end{array}\right]\,.
\end{equation}
In order to highlight nilpotency, we  introduce the \addedinrevision{diagonal cancellation matrices $I_l$ (resp. $I_r$) such that $(I_l)_{ii}=1$ if $1\le i\le N-1$ (resp. $(I_r)_{ii}=1$  if $N\le i \le 2N-2$) and $0$ otherwise,} as well as four $(2N-2)\times (2N-2)$ submatrices of the operator ${\mathscr T}_{LR}$:
$$ {\mathscr M}_l := I_l^T\, {\mathscr T}_{LR}\, I_l\ \text{(top left)}\,,\ {\mathscr A}_l := I_l^T\, {\mathscr T}_{LR}\, I_r\  \text{(top right)}$$
$$ {\mathscr A}_r := I_r^T\, {\mathscr T}_{LR}\, I_l\   \text{(bottom left)}\,,\ {\mathscr M}_r := I_r^T\, {\mathscr T}_{LR}\, I_r \ \text{(bottom right )}$$ 
so that we have $I-{\mathscr T}_{LR}=I-{\mathscr M}_l-{\mathscr A}_l-{\mathscr M}_r-{\mathscr A}_r$. It is easy to check that we have the following cancellation relations:
\begin{equation}
\boxed{
\begin{array}{lll}
{\mathscr M}_r^{N-1}={\mathscr M}_l^{N-1}=0;&{\mathscr M}_l\,{\mathscr M}_r={\mathscr M}_r\,{\mathscr M}_l=0;&{\mathscr A}_l^2={\mathscr A}_r^2=0\\
{\mathscr A}_l\,{\mathscr M}_l={\mathscr A}_r\,{\mathscr M}_r=0;&{\mathscr M}_l\,{\mathscr A}_r={\mathscr M}_r\,{\mathscr A}_l=0\,.&\ \ \ \\
\end{array}
}\label{eq:tabref}
\end{equation}
It is worth noticing that these relations come from
the structure of the matrices and do not depend on the value of the
entries. Let us introduce the following operators that will play a crucial role in the sequel:
\begin{equation}
	\label{eq:CrCl}
	\boxed{	C_r := (I-{\mathscr M}_r)^{-1}\,{\mathscr A}_r =  (\sum_{i=0}^{N-2} {\mathscr M}_r^i)\,{\mathscr A}_r \text{ \ and \ } C_l := (I-{\mathscr M}_l)^{-1}\,{\mathscr A}_l = (\sum_{i=0}^{N-2} {\mathscr M}_l^i)\,{\mathscr A}_l\,.}
\end{equation}
Note that using cancellation relations~\eqref{eq:tabref}, we have 
\begin{equation}
	\label{eq:CR2CL2zero}
C_r^2 = C_l^2 = 0 \reviewerone{\,; {\mathscr M}_r\,C_l = {\mathscr M}_l\,C_r = 0 \text{ and } {\mathscr A}_r\,C_r={\mathscr A}_l\,C_l=0}\,.
\end{equation}

\subsection{Jacobi, Gauss-Seidel and Symmetric Gauss-Seidel for the subdomain wise numbering}
	\label{sec:algosSW}
	
With the notations introduced above, the substructured operator~\eqref{eq:substructSW} reads:
\begin{equation}
	\label{eq:substructSWAMnotations}
 (I-{\mathscr T}_{SW})=
\left[ 
\begin{array}{cccccccc}
   I&-{\mathscr A}_{r,2N-2,1} &-{\mathscr M}_{r,2N-2,2N-3} & 0 &\phantom{{-\mathscr T}}&\phantom{{-\mathscr T}} &  &  \\
   -{\mathscr A}_{l,1,2N-2}& I & 0 & 0 &\ddots &\ &\ &\  \\
   0 & 0 & I &-{\mathscr A}_{r,2N-3,2} &-{\mathscr M}_{r,2N-3,2N-4} & 0 &\ &\  \\
   \ &{\mathscr M}_{l,2,1} &{\mathscr A}_{l,2,2N-3}& I & 0 & 0 &\ &\  \\
   \ & \ & 0 & 0 &\ddots &\ddots &\ddots &\  \\
   \ &\ &\ &-{\mathscr M}_{l,3,2}  &\ddots &\ddots &\ddots &\ddots  \\
   \ &\ &\ &\ &\ddots &\ddots &\ddots &\ddots  \\
   \ &\ &\ &\ &\ &\ddots\ &\ddots\ &I   
\end{array}
\right]
\,.
\end{equation} 	
When considering the decomposition of matrix $I-{\mathscr T}_{SW}$ into its diagonal, lower and upper parts as $I-{\mathscr T}_{SW}=L+D+U$, its upper part $U$ involves only entries of the matrices ${\mathscr A}_r$ and ${\mathscr M}_r$ and similarly, its lower part $L$ involves only entries of the matrices ${\mathscr A}_l$ and ${\mathscr M}_l$ and the diagonal part $D$ is simply the identity. {\em Expressed now in the LR numbering}, it is then easy to check that we have $D=I$, $L=-{\mathscr M}_l-{\mathscr A}_l$ and $U=-{\mathscr M}_r-{\mathscr A}_r$. Thus we introduce in a classical way several preconditioners whose formulas are given in the LR numbering:
\begin{itemize}
\item the Jacobi preconditioner $D^{-1}$: 
\[
M^{-1}_{Jacobi} := I \,,
\]
\item the left Gauss-Seidel (GS)  preconditioner $(L+D)^{-1}$:
\[
	M^{-1}_{GS} := (I-{\mathscr M}_l- {\mathscr A}_l)^{-1}\,,
\]
\item the Symmetric Gauss-Seidel (SGS)  preconditioner $[(L+D) D^{-1} (D+U)]^{-1} = (D+U)^{-1}\,D\,(L+D)^{-1}$:
\[
M^{-1}_{SGS}  := (I-{\mathscr M}_r- {\mathscr A}_r)^{-1}\,(I-{\mathscr M}_l- {\mathscr A}_l)^{-1} \,.
\]
\end{itemize}
The Jacobi preconditioner amounts to Desprès-Lions algorithm and the two other preconditioners had been introduced and studied in~\cite{Nataf:1997:CRO,vion:2018:improved}. In~\cite{Nataf:1997:CRO}, the $GS$ preconditioner was referred to as FDA and the $SGS$ preconditioner as the Double Sweep (DS) algorithm.

When used as solvers, their convergence depends on the spectral radius of the error propagation operators defined as:
\begin{madefinition}
	\label{def:errorpropagation}
For a preconditioner $M^{-1}$ of a linear operator $\reviewerone{{\mathscr T}}$, we denote by $R$ its error propagation operator as: $R:=M^{-1}(M-\reviewerone{{\mathscr T}})$. 
\end{madefinition}
It was proved in~\cite{Nataf:1997:CRO} that the spectral radius of their error propagation operators denoted by $\rho(R)$ satisfy the following bounds: $\rho(R_{Jacobi})\le\rho^{1/N}$, $\rho(R_{GS})\le\rho^{2/N}$ and $\rho(R_{SGS})\le\rho$  where $N$ is the number of subdomains and $\rho=\rho(C_r)\,\rho(C_l)$. The proof was written with the help of the formal language theory, see e.g.~\cite{lothaire:1997:combinatorics}. The $SGS$ algorithm can be studied more classically since we have:
\begin{proposition}
	\label{th:cvopSGS}
The formula for $R_{SGS}$ is:
\[
	R_{SGS} = (I\, +C_r)\,C_l\,({\mathscr M}_r+{\mathscr A}_r)\,.
\]
As for the power of $R_{SGS}$, we have the following formula:
\[ 
	  R_{SGS}^n = (I\, +C_r)\,C_l\,(C_r\,C_l)^{n-1}\,({\mathscr M}_r+{\mathscr A}_r)\,.
\]
\end{proposition} 
\begin{proof}
By definition, we have:
\[
	  R_{SGS} = (Id-{\mathscr A}_r-{\mathscr M}_r)^{-1} (Id-{\mathscr M}_l-{\mathscr A}_l)^{-1}({\mathscr M}_l+{\mathscr A}_l) ({\mathscr A}_r+{\mathscr M}_r)\,.
\]
\reviewerone{ Using cancellations relations~\eqref{eq:tabref}, we have:
\[
	Id-{\mathscr A}_l-{\mathscr M}_l = (Id-{\mathscr A}_l)\,(Id-{\mathscr M}_l) 
	\text{ and } (Id-{\mathscr A}_l)^{-1} = Id+{\mathscr A}_l\,,
\]
so that 
\[
	(Id-{\mathscr A}_l-{\mathscr M}_l)^{-1} = (Id-{\mathscr M}_l)^{-1}\,(Id+{\mathscr A}_l)\,.
\]
Similarly, we have:
\[
	(Id-{\mathscr A}_r-{\mathscr M}_r)^{-1} = (Id-{\mathscr M}_r)^{-1}\,(Id+{\mathscr A}_r)\,.
\]
Thus, the expression for $R_{SGS}$ simplifies as:
\[
	R_{SGS} = (Id-{\mathscr M}_r)^{-1}\,(Id+{\mathscr A}_r)(Id-{\mathscr M}_l)^{-1}\,(Id+{\mathscr A}_l)({\mathscr M}_l+{\mathscr A}_l) ({\mathscr A}_r+{\mathscr M}_r)\,.
\]
Using ${\mathscr A}_l\,({\mathscr M}_l+{\mathscr A}_l)=0$ and then ${\mathscr M}_l\,({\mathscr A}_r+{\mathscr M}_r)=0$, we get 
\[
	R_{SGS} = (Id-{\mathscr M}_r)^{-1}\,(Id+{\mathscr A}_r)(Id-{\mathscr M}_l)^{-1}\,{\mathscr A}_l ({\mathscr A}_r+{\mathscr M}_r)
	= ((Id-{\mathscr M}_r)^{-1} + C_r)\,C_l\,({\mathscr A}_r+{\mathscr M}_r)\,.
\]
}
Note that since ${\mathscr M}_r\,C_l=0$ and thus $(Id-{\mathscr M}_r)^{-1}\, C_l=C_l$ \reviewerone{ we have that the expression for $R_{SGS}$ simplifies as:
\[
	R_{SGS} =(Id + C_r)\,C_l\,({\mathscr A}_r+{\mathscr M}_r)\,.
\]
}\\
As for the powers of $R_{SGS}$, we first note that:
	\[
	  R_{SGS}^n = (I\, +C_r)\,C_l\,\left[({\mathscr M}_r+{\mathscr A}_r) (I\, +C_r)\,C_l\right]^{n-1} ({\mathscr M}_r+{\mathscr A}_r)  \,.
	\]
\reviewerone{
Using cancellation relations~\eqref{eq:tabref} and~\eqref{eq:CR2CL2zero} to simplify the middle term, we get:
	\[
	\begin{array}{rcl}
	({\mathscr M}_r+{\mathscr A}_r) (I\, +C_r)\,C_l 
	 &=& ({\mathscr M}_r+{\mathscr A}_r) (I\, +C_r)\,C_l = ({\mathscr M}_r+{\mathscr A}_r + {\mathscr M}_r\,C_r)\,C_l    \\
	 &=& \left[ (\sum_{i=1}^{N-2}  {\mathscr M}_r^i )\,{\mathscr A}_r + {\mathscr M}_r + {\mathscr A}_r \right] C_l = \left[ C_r + {\mathscr M}_r  \right] C_l  \\
	 &=& C_r \,C_l\,.
	\end{array}
	\]
}
\end{proof}

\subsection{Block Jacobi, Gauss-Seidel and Symmetric Gauss-Seidel for the left-right numbering}
	\label{sec:algosLR}

Considering system~\eqref{eq:substructLR} as a block two by two matrix, the matrix $I-{\mathscr T}_{LR}$ can be decomposed into its block diagonal, lower and upper parts as $L+D+U$ where $D=I-{\mathscr M}_r-{\mathscr M}_l$, $L=-{\mathscr A}_r$ and $U=-{\mathscr A}_l$. Thus we introduce in a classical way several preconditioners:
\begin{itemize}
\item the Block Jacobi  preconditioner (BJ) $D^{-1}$ :
\begin{equation}
	\label{eq:precondBlockLRBJ}
\begin{array}{rcl}
	M^{-1}_{BJ} &:=& (I-{\mathscr M}_r-{\mathscr M}_l)^{-1}\,, \\ 
\end{array}	
\end{equation}
\item the Block Gauss-Seidel  (BGS)  preconditioner $(L+D)^{-1}$:
\begin{equation}
	\label{eq:precondBlockLRBGS}
\begin{array}{rcl}
	M^{-1}_{BGS} &:=& (I-{\mathscr M}_r-{\mathscr M}_l- {\mathscr A}_r)^{-1}\,, \\ 
\end{array}	
\end{equation}
\item the Block Symmetric Gauss-Seidel  (BSGS)  preconditioner $[(L+D) D^{-1} (D+U)]^{-1}= (D+U)^{-1}\,D\,(L+D)^{-1}$:
\begin{equation}
	\label{eq:precondBlockLRBSGS}
\begin{array}{rcl}
	M^{-1}_{BSGS}  &:=& (I-{\mathscr M}_r-{\mathscr M}_l- {\mathscr A}_l)^{-1}\,(I-{\mathscr M}_r-{\mathscr M}_l)\,(I-{\mathscr M}_r-{\mathscr M}_l- {\mathscr A}_r)^{-1} \,.
\end{array}	
\end{equation}
\end{itemize}

\begin{remarque}
The BJ preconditioner was introduced in~\cite{vion:2014:double} for non overlapping decompositions. Note that a diagonal block solve consists in two independent forward substitutions on the left and right interface unknowns. So the method was called Double Sweep although it is different from the one previously introduced in~\cite{Nataf:1997:CRO} and named as well Double Sweep. We hope to have clarified this possible confusion. \reviewertwo{The BGS preconditioner was introduced in~\cite{Stolk:2013:rapidly,Stolk:2017:improved} for overlapping decompositions. Note that it corresponds to the preconditioner $(I-\mathscr{T}^{EABC}_{LR,GS})^{-1}$ introduced in \S~\ref{sec:threesditerative}.}
\end{remarque}

\smallskip

To the best of our knowledge, the BSGS algorithm has not been previously introduced in the literature. Neither its convergence rate nor that of BJ and BGS has been studied before. The emergence of the BSGS algorithm stems from our approach to introduce the BJ algorithm, which differs from the one presented in~\cite{vion:2014:double} and has been previously outlined in~\S~\ref{sec:stolkviongeuzaine}.

We have the following propositions.
\begin{proposition}
	\label{th:cvopBJ}
The error propagation operator $R_{BJ}$ has the following expression:
\[
R_{BJ} = C_r + C_l\,,
\]
and for $n$ even we have:	
\[
R_{BJ}^n = (C_r\,C_l)^{n/2}   + (C_l\,C_r)^{n/2} \,.
\]
\end{proposition}
\begin{proof}
Thanks to relations~\eqref{eq:tabref}, we have $M_{BJ}^{-1}=Id+\sum_{i=1}^{N-2}{\mathscr M}_l^i+\sum_{i=1}^{N-2}{\mathscr M}_r^i$ and:
\[
	R_{BJ} = (Id+\sum_{i=1}^{N-2}{\mathscr M}_l^i+\sum_{i=1}^{N-2}{\mathscr M}_r^i)({\mathscr A}_l+{\mathscr A}_r) = \sum_{i=0}^{N-2} {\mathscr M}_r^i\,{\mathscr A}_r + \sum_{i=0}^{N-2} {\mathscr M}_l^i\,{\mathscr A}_l = C_r + C_l\,.
\]
Next using cancellation relation~\eqref{eq:CR2CL2zero}, the formula for $R_{BJ}^n$ can easily be proved by induction.
\end{proof}

\begin{proposition}
	\label{th:cvopBGS}
The error propagation operator $R_{BGS}$ has the following expression:
\[
R_{BGS} = (I+ C_r)\,C_l\,\,,
\]
and for any integer $n$ we have:	
\[
R_{BGS}^n = (I+C_r)\,C_l (C_r\,C_l)^{n-1} \,.
\]	
\end{proposition}
\begin{proof}
\reviewerone{
We first prove the following equalities:
\begin{lemma}
\begin{equation}
	\label{eq:cancel2}
\begin{array}{rcl} 
(I-{\mathscr M}_r-{\mathscr M}_l)^{-1}  &=& (I-{\mathscr M}_r)^{-1}\,(I-{\mathscr M}_l)^{-1} = (I-{\mathscr M}_l)^{-1}\,(I-{\mathscr M}_r)^{-1} \\
(I-{\mathscr M}_r-{\mathscr M}_l-{\mathscr A}_r)^{-1}&=&(I-{\mathscr M}_r)^{-1}\,(I+{\mathscr A}_r)\,(I-{\mathscr M}_l)^{-1} \\
(I-{\mathscr M}_r-{\mathscr M}_l-{\mathscr A}_l)^{-1}&=&(I-{\mathscr M}_l)^{-1}\,(I+{\mathscr A}_l)\,(I-{\mathscr M}_r)^{-1}\,,
\end{array}
\end{equation}
\end{lemma}
\begin{proof}
Since ${\mathscr M}_r\,{\mathscr M}_l={\mathscr M}_l\,{\mathscr M}_r=0$, we have 
\[
  (I-{\mathscr M}_r-{\mathscr M}_l) = (I-{\mathscr M}_r)(I-{\mathscr M}_l) = (I-{\mathscr M}_l)(I-{\mathscr M}_r)\,,
\]
and the first equality of the lemma follows by taking the inverse of the above equation.\\
Using cancellation relations~\eqref{eq:tabref}, we have:
\[
(I-{\mathscr M}_r-{\mathscr M}_l-{\mathscr A}_r) = (I-{\mathscr M}_l)(I-{\mathscr M}_r-{\mathscr A}_r)
= (I-{\mathscr M}_l)(I-{\mathscr A}_r)(I-{\mathscr M}_r)
\]
Taking the inverse of this relation and using $(I-{\mathscr A}_r)^{-1}=I+{\mathscr A}_r$ proves the second equation of the lemma. The proof of the third line is similar. 
\end{proof}

We will now resume the proof of the proposition:
}
\[
\begin{array}{rcl} 
R_{BGS}  &=& M^{-1}_{BGS}\,{\mathscr A}_l = (I-{\mathscr M}_r)^{-1}\,(I+{\mathscr A}_r)\,(I-{\mathscr M}_l)^{-1}\,{\mathscr A}_l\\
 &=& (I-{\mathscr M}_r)^{-1}\,(I-{\mathscr M}_l)^{-1}\,{\mathscr A}_l +  (I-{\mathscr M}_r)^{-1}\,{\mathscr A}_r\,(I-{\mathscr M}_l)^{-1}\,{\mathscr A}_l \\
&=& (I-{\mathscr M}_l)^{-1}\,(I-{\mathscr M}_r)^{-1}\,{\mathscr A}_l +  C_r\,C_l = C_l +  C_r\,C_l \,.
\end{array}
\]
As for the $n$-th power of $R_{BGS}$, the formula can be proved by induction using the fact that:
\[
(I+ C_r)\,C_l\, (I+ C_r)\,C_l\, = (I+ C_r)\,C_l\,  C_r\,C_l\,.
\]
\end{proof}

\begin{proposition}
	\label{th:cvopBSGS}
The error propagation operator $R_{BSGS}$ has the following expression:
\[
R_{BSGS} = (I+ C_l)\,C_r\, C_l\,\,,
\]
and for any integer $n$ we have:	
\[
R_{BSGS}^n = (I+C_l)\, (C_r\,C_l)^{n} \,.
\]
\end{proposition}
\begin{proof}
Using cancellation relations~\eqref{eq:tabref} and \eqref{eq:cancel2}, we have:
\[
\begin{array}{rcl} 
	R_{BSGS} &=& M^{-1}_{BGS} L\,D^{-1}\,U = M^{-1}_{BGS} {\mathscr A}_r\,(I-{\mathscr M}_r)^{-1}\,(I-{\mathscr M}_l)^{-1}\,{\mathscr A}_l = M^{-1}_{BGS}\,{\mathscr A}_r\,C_l\\
	  &=& (I-{\mathscr M}_l)^{-1}\,(I+{\mathscr A}_l)\,(I-{\mathscr M}_r)^{-1} (I-{\mathscr M}_r)(I-{\mathscr M}_l)(I-{\mathscr M}_r)^{-1}\,(I+{\mathscr A}_r)\,(I-{\mathscr M}_l)^{-1}\,{\mathscr A}_r\,C_l\\
	&=& (I-{\mathscr M}_l)^{-1}\,(I+{\mathscr A}_l-{\mathscr M}_l)(I-{\mathscr M}_r)^{-1}\,(I+{\mathscr A}_r)\,(I-{\mathscr M}_l)^{-1}\,{\mathscr A}_r\,C_l\\
	&=& (I+C_l)\,(I-{\mathscr M}_r)^{-1}\,(I+{\mathscr A}_r)\,{\mathscr A}_r\,C_l = (I+C_l)\,C_r\,C_l\,.
\end{array}
\]
As for the last formula of the proposition, it comes from the nullity of both $C_r^2$ and $C_l^2$.
\end{proof}

\begin{remark}
It is clear from the above results that the operators $C_r$ and $C_l$ (see eq.~\eqref{eq:CrCl}) are a key measure of the efficiency of these algorithms. If the left interface conditions are exact absorbing conditions, operator ${\mathscr A}_r=0$, see equation~\eqref{eq:AzeroEABC} and thus $C_r=0$ as well. Similarly if right interface conditions are exact absorbing conditions, operator $C_l=0$. More generally, the norms of $C_r$ and $C_l$ are proportional to that of ${\mathscr A}_r$ and ${\mathscr A}_l$. Thus as expected the more absorbing the interface conditions are, the better the convergence is. Another parameter is the number of subdomains since for given interface conditions, as the number of subdomains increases, the norm of $C_r$ and $C_l$ will grow. It echoes what was noticed in~\cite{Stolk:2013:rapidly} at the end of \S~4: 
\begin{center}
{\em "With the Robin transmission conditions the iteration numbers grow roughly linearly in $N_x$, or as $N^{1/2}$ in 2-D."}
\end{center}

In this respect, note that $R_{BSGS}$ satisfies
\[
	(I+ C_l)^{-1}\,R_{BSGS}\,(I+ C_l) = C_r\, C_l\,
\]
so that from all algorithms listed in Table~\ref{tab:algocvrate} it has the most favorable amplification error operator iff the norm of $C_r\,C_l$ is smaller than one. In a recent article~\cite{Gong:2021:Convergence}, these norms have been estimated in the context of the study of the convergence rate of the Jacobi preconditioner. 
\end{remark}

\subsection{Theoretical comparison between the SGS and BJ preconditioners}

We devote a paragraph to this comparison since these algorithms have been misleadingly coined with the same name, i.e. Double Sweep algorithm, in the articles in which they were introduced, see~\cite{Nataf:1997:CRO,vion:2018:improved} for SGS and~\cite{vion:2014:double} for BJ.  
From Propositions~\ref{th:cvopSGS} and ~\ref{th:cvopBJ}, we have for $n$ even:
\[
R_{BJ}^n = (C_r\,C_l)^{n/2} + (C_l\,C_r)^{n/2}\,,
\] 
and for any $n$
\[
R_{SGS}^n =  (I+C_r)\,C_l\, (C_r\,C_l)^{n-1}\,({\mathscr M}_r+{\mathscr A}_r) \,.
\]
Let us denote by 
\begin{equation}
	\label{eq:rhodef}
\rho \text{ the spectral radius of } C_r\,C_l\,,
\end{equation}
which is the same as that of $C_l\,C_r$. We have for any operator norm $\|\cdot\|$ on the matrices the following estimates:
\begin{equation}
	\| R_{BJ}^n \| \le  \| (C_r\,C_l)^{n/2} \| + \| (C_l\,C_r)^{n/2} \| \le \| (C_r\,C_l)^{n/2} \| + \| C_l \| \| (C_r\,C_l)^{n/2-1} \| \| C_l \|   \,,
\end{equation}
so that taking the $n$-th square root, we get the following estimate for the spectral radius of $R_{BJ}$:
\begin{equation}
	\label{eq:cvEstimateGlobalBJb}
	\rho(R_{BJ}) \le \rho^{1/2}\,.
\end{equation}
As for the SGS algorithm, we get:
\begin{equation}
	\| R_{SGS}^n \| \le \| I+C_r\|\, \|C_l \|\,	\| {\mathscr M}_r+{\mathscr A}_r \|\,\, \|(C_r\,C_l)^{n-1}\| \,\,.	
\end{equation}
We thus get the following estimate for the spectral radius of $R_{SGS}$: 
\begin{equation}
	\label{eq:cvEstimateGlobalSGSb}
	\rho(R_{SGS}) \le  \rho\,.	
\end{equation}
 Thus when $\rho < 1$, an advantage of SGS over BJ is the square root of $\rho$ in the ratio of the convergence rates. This is coherent with the factor two in the iteration counts in favour of SGS  observed in most tables of the numerical section~\ref{sec:resnum} when PMLs are interface conditions. \reviewertwo{Note that the preconditioner BGS has the same convergence properties as the preconditioner SGS.}
 \begin{table}[htp]
 \small
\begin{tabular}{|l|c|l|l|l|l|}\hline
	\bfseries Ref. & \bfseries Abbr. & \bfseries Linear Algebra & \bfseries Definition & \bfseries Ampl. Error & \bfseries Radius
\begin{filecontents*}{algocvratetalk.csv}
noms , numbering , linalg ,      definition ,     cvrate ,     asymptoticcvrate  
$\begin{array}{l} \text{Despr\'es (1991)} \end{array}$  , Jacobi  ,  Jacobi-SW	,  $I$ ,  ${\mathscr M}_r+{\mathscr M}_l+{\mathscr A}_r+{\mathscr A}_l$  ,   $\ \rho^{1/N}$
$\begin{array}{l} \text{Nataf-Nier (1997)} \end{array}$ , GS  ,  Gauss-Seidel-SW	,  $(I-{\mathscr M}_l-{\mathscr A}_l)^{-1}$ ,  $(I+C_l)({\mathscr M}_r+{\mathscr A}_r)$  ,   $\ \rho^{2/N}$
$\left\{\begin{array}{l} \text{Nataf-Nier (1997)} \\ \text{Vion-Geuzaine(2018)} \end{array}\right.$  , SGS  ,  Symm. Gauss-Seidel-SW , $\begin{array}{l} (I-{\mathscr M}_r-{\mathscr A}_r)^{-1} \\ \ \ \  \times (I-{\mathscr M}_l-{\mathscr A}_l)^{-1} \end{array}$ , $(I\, +C_r)\,C_l\,({\mathscr M}_r+{\mathscr A}_r)$ , $\ \rho$
$\begin{array}{l}  \text{Vion-Geuzaine (2014)}\end{array}$   , BJ  ,  Block Jacobi-LR	,  $(I-{\mathscr M}_r-{\mathscr M}_l)^{-1}$ ,  $C_r+C_l$  ,   $\rho^{1/2}$
$\begin{array}{l} \text{\reviewertwo{Stolk (2013)}} \end{array}$     , BGS  ,  Block GS-LR  , $(I-{\mathscr M}_r-{\mathscr M}_l-{\mathscr A}_r)^{-1}$ , $(I+ C_r)\,C_l$ , $\ \rho$
$\begin{array}{l} \text{Here (2022)} \end{array}$    , BSGS  ,  Symm. Block  GS-LR , $\begin{array}{r} (I-{\mathscr M}_r-{\mathscr M}_l-{\mathscr A}_l)^{-1}\\ \times (I-{\mathscr M}_r-{\mathscr M}_l)\\ \times(I-{\mathscr M}_r-{\mathscr M}_l-{\mathscr A}_r)^{-1} \end{array}$ , $(I+ C_l)\,C_r\,C_l$ , $\ \rho$

\end{filecontents*}
\csvreader[head to column names]{algocvratetalk.csv}{}
{\\[.95em] \noms & \numbering & \linalg & \definition & \cvrate & \asymptoticcvrate}\\
\hline
\end{tabular}
\caption{ Algorithms and their convergence properties. SW means subdomain wise numbering and LR means left right numbering }
 \label{tab:algocvrate}
 \end{table}


\reviewertwo{
\subsection{Cost evaluation of the preconditioners}

	\label{sec:costevaluation}
	
	\reviewertwo{Before considering the cost of the preconditioners, note that they are all used to solve the substructured system~\eqref{global_substructuring} and the residual has to be computed at every iteration. The matrix-vector product with the operator ${\mathscr T}$ amounts to one subdomain solve in parallel on all subdomains. This cost has to be added to that of the preconditioners that we estimate now.\\}

For instance, the SGS preconditioner necessitates two solves per subdomain, as outlined in the volumic formulation in equations~\eqref{leftsweep_algo_helmholtz}-\eqref{rightsweep_algo_helmholtz}. It can also be derived at the substructured level. Indeed, consider the application of $(I-{\mathscr M}_l-{\mathscr A}_l)^{-1}$. It amounts to solve the following linear system for some right hand side $G$:
\begin{equation}
	\label{eq:substructSGS1solve}
\left[ \begin{array}{c;{2pt/2pt}c}
\begin{array}{cccc}
   I& \ & \ & 0 \\
   -{\mathscr T}_{(3,l)(2,l)}&\ddots& \ &\ \\
   \ &\ddots&\ddots&\ \\
   0&\ &-{\mathscr T}_{(N,l)(N-1,l)}&I 
   \end{array}&
\begin{array}{cccc}
  0&\ &\ &-{\mathscr T}_{(2,l)(1,r)}\\
   \ &\ &\iddots &\ \\
   \ &\iddots &\ &\ \\
   -{\mathscr T}_{(N,l)(N-1,r)}&\ &\ &0 
   \end{array}\\ \hdashline[2pt/2pt]
\begin{array}{cccc}
  &\ &\ &\\
   \ &\ & &\ \\
   \ & &0\ &\ \\
   &\ &\ & 
   \end{array} &  
   \begin{array}{cccc}
   I& &\ \ \ &\phantom{-{\mathscr T}_{(N,l)(N-1,r)}}  \\[.5em]
   &\ddots&\ &\ \\[.5em]
   \ &\ &\ddots&\ \ \\[.5em]
   \phantom{-{\mathscr T}_{(N,l)(N-1,r)}}  \ &\ \ &\ \  &I 
   \end{array}
\end{array}\right]\ 
\left[\begin{array}{c} h_{2,l}\\ \vdots \\ \vdots \\ h_{N,l}\\
h_{N-1,r}\\ \vdots\\ \vdots \\ h_{1,r}
\end{array}\right]\,
= 
\left[\begin{array}{c} g_{2,l}\\ \vdots \\ \vdots \\ g_{N,l}\\
g_{N-1,r}\\ \vdots\\ \vdots \\ g_{1,r}
\end{array}\right]\,.
\end{equation}
From the bottom rows, we have directly $h_{i,r} = g_{i,r}$ for all $1\le i\le N-1$. Then the forward sweep on the upper rows reads:
\[
 h_{i,l} = {\mathscr T}_{(i,l)(i-1,l)}(h_{i-1,l}) + {\mathscr T}_{(i,l)(i-1,r)}(g_{i-1,r})\,.
\]
The above quantity can be evaluated with only one solve in subdomain~$i-1$ since by definitions~\eqref{eq:Tentries} and the linearity of operator $S_{i-1}$, we have:
\[ 
h_{i,l} = {\mathscr B}_{i,l}(S_{i-1}(h_{i-1,l},0,0)) + {\mathscr B}_{i,l}(S_{i-1}(0,g_{i-1,r},0))
 = {\mathscr B}_{i,l}(S_{i-1}(h_{i-1,l},g_{i-1,r},0))\,.
\]
With the same kind of reasoning, we see that the cost of the preconditioners is as follows:\\ Jacobi: 0 solve \ ; \ GS  N solves\ ;  SGS  2N solves\ ; \ BJ 2N solves \ ; \ BGS   2N solves\ ; \ BSGS 	4N solves.\\
} 
Note that at the expense of doubling the number of cores used in the application of the BJ preconditioner, the difference in the convergence rate estimate w.r.t. to SGS and BGS preconditioners is nullified. Indeed, during one iteration of the SGS or BGS algorithms, only one subdomain is active at a time in the order $1\rightarrow 2 \rightarrow \ldots \rightarrow N$ and then in the reverse order $N \rightarrow N-1 \rightarrow \ldots \rightarrow 1$. Whereas during one iteration of BJ, two subdomains are active at a time in the order $1\rightarrow 2 \rightarrow \ldots \rightarrow N$ for one subdomain and concurrently in the reverse order $N \rightarrow N-1 \rightarrow \ldots \rightarrow 1$ for the other one. Thus, the elapsed time of one iteration of SGS or BGS on one core is twice as much as one iteration of BJ using two cores. This might compensate for the higher iteration counts of BJ compared to SGS or BGS when the number of cores is not limited.\\

\subsection{Volumic preconditioners} 
\label{sub:volumicalgo}

The substructured algorithms demand exact local solves and data structures adapted to substructured algorithms. In order to allow for the use of fast approximate solvers in subdomains and of volumic data structures, an extension to the volumic case is defined below. For the SGS algorithm, the volumic formulation is simply given by equations~\eqref{leftsweep_algo_helmholtz}-\eqref{rightsweep_algo_helmholtz} followed by the use of a partition of unity in order to iterate on functions which are uniquely defined in the overlaps. For this, we first introduce for $1\le i \le N$ partition of unity functions $\chi_i:\Omega_i \rightarrow \mathbb{R}^+$ and for $v_i:\Omega_i \rightarrow \mathbb{C}$, $E_i(v_i)$ denotes the extension by zero of $v_i$  to $\Omega$ so that for any function $v:\Omega \rightarrow \mathbb{C}$, we have:
\[
	v = \sum_{i=1}^N E_i(\chi_i\,v_{|\Omega_i})\,.
\]
We introduce the operator $MC$ (MakeCoherent) that maps a collection of local functions $(v_i)_{1\le i\le N}$ to a global function defined as follows:  
\begin{equation}
  \label{eq:makecoherent}
  MC((v_i)_{1\le i\le N}) := \sum_{i=1}^N E_i( \chi_i\,v_i)\,.
\end{equation}
Let $f$ be a source term and $h$ be interface sources, it is then natural to introduce the following linear parallel reconstruction algorithm $\text{Vol}$:
\begin{equation}
  \label{eq:operatorVol}
  \text{Vol}(h,f) := MC((S_i(h_{i,l},h_{i,r},f))_{1\le i \le N})\,.
\end{equation}
Let $u$ be the solution to the original problem~\eqref{helmholtz_equation}, we have:
\begin{equation}
  \label{eq:voltosub}
  u = \text{Vol}((I-{\mathscr T})^{-1}G(f),f)\,.
\end{equation}
Since a substructured preconditioner $M^{-1}$ may be seen as an approximate inverse to $(I-{\mathscr T})$, it is natural to define the related volumic precondtioner as:
\begin{madefinition}
	\label{def:volprec}
Let $M^{-1}$ be a substructured preconditioner for problem~\eqref{global_substructuring}, the related volumic preconditioner to problem~\eqref{helmholtz_equation} $M_{vol}^{-1}$ is defined as follows:
\begin{equation}
  \label{eq:volprec}
  M_{vol}^{-1}(f) := \text{Vol}(M^{-1}G(f),f) = \text{Vol}(M^{-1}G(f),0) + \text{Vol}(0,f)\,.
\end{equation}
\end{madefinition}
Note that if the surface variant of the preconditioner is exact as it happens for some if the interface conditions are exact absorbing boundary conditions, this property will be inherited by its volumic counterpart. \\

As an example, we detail below after some rewriting the action of the volumic preconditioner associated to the surface variant of the BJ preconditioner.  

\begin{minipage}{\textwidth}
 \begin{center}
 	{\bf Volumic BJ Preconditioner}
 \end{center}	
\begin{multicols}{2}
{\bf Left to right sweep}\\
 Subdomain $1$\\
$w_1 \leftarrow S_1(0,0,f)$\\
 Subdomain $2$\\
$G_{2,l}(f) \leftarrow \mathcal{B}_{2,l}(w_1)$\\
$h_{2,l} \leftarrow G_{2,l}(f)$\\
$v_{2,l} \leftarrow S_2(h_{2,l},0,0)$\\
$w_2 \leftarrow S_2(0,0,f)$\\
 Subdomain $3$\\
$G_{3,l}(f) \leftarrow \mathcal{B}_{3,l}(w_2)$\\
$h_{3,l} \leftarrow G_{3,l}(f)+\mathcal{B}_{3,l}(v_{2,l})$\\
$v_{3,l} \leftarrow S_3(h_{3,l},0,0)$\\
$w_3 \leftarrow S_3(0,0,f)$\\
 Subdomain $4$\\
$G_{4,l}(f) \leftarrow \mathcal{B}_{4,l}(w_3)$\\
$h_{4,l} \leftarrow G_{4,l}(f)+\mathcal{B}_{4,l}(v_{3,l})$\\
\hphantom{$h_{4,l} \leftarrow$} $\vdots$ \\
 Subdomain $N$\\
$w_N \leftarrow S_N(0,0,f)$\\
$v_{N,l} \leftarrow S_N(h_{N,l},0,0)$\\

	\columnbreak
	
{\bf Right to left sweep}\\
 Subdomain $N$\\
$w_N \leftarrow S_N(0,0,f)$\\
 Subdomain $N-1$\\
$G_{N-1,r}(f) \leftarrow \mathcal{B}_{N-1,r}(w_N)$\\
$h_{N-1,r} \leftarrow G_{N-1,r}(f)$\\
$v_{N-1,r} \leftarrow S_{N-1}(0,h_{N-1,r},0)$\\
$w_{N-1} \leftarrow S_{N-1}(0,0,f)$\\
 Subdomain $N-2$\\
$G_{N-2,r}(f) \leftarrow \mathcal{B}_{N-2,r}(w_{N-1})$\\
$h_{N-2,r} \leftarrow G_{N-2,r}(f)+\mathcal{B}_{N-2,r}(v_{N-1,r})$\\
$v_{N-2,r} \leftarrow S_{N-2}(0,h_{N-2,r},0)$\\
$w_{N-2} \leftarrow S_{N-2}(0,0,f)$\\
 Subdomain $N-3$\\
$G_{N-3,r}(f) \leftarrow \mathcal{B}_{N-3,r}(w_{N-2})$\\
$h_{N-3,r} \leftarrow G_{N-3,r}(f)+\mathcal{B}_{N-3,r}(v_{N-2,r})$\\
\hphantom{$h_{4,l} \leftarrow$} $\vdots$ \\
 Subdomain $1$\\
$w_1 \leftarrow S_1(0,0,f)$\\
$v_{1,r} \leftarrow S_1(h_{1,r},0,0)$\\
	\end{multicols}
\end{minipage}	
Then we perform local summations:
\[
	u_i := w_i+v_{i,l}+v_{i,r}, \ \ \ \ 2 \le i \le N-1
\]
and
\[
  u_1:= w_1 + v_{1,r}\ \ \text{ and } u_N := w_N + v_{N,l}\,.
\]
Finally the action of the preconditioner $M^{-1}_{BJ,vol}$ reads:
\begin{equation}
\label{eq:volpreconditioner}
  M^{-1}_{BJ,vol}(f) := \sum_{i=1}^N E_i(\chi_i (u_i))\,. 
\end{equation}

%
%
%
%

\section{Numerical results}
	\label{sec:resnum}

In this section, we present numerical results for the volumic variants of the preconditioners listed in Table~\ref{tab:algocvrate}, used with the GMRES algorithm~\cite{Saad:1986:GGM}. We consider two different relative tolerances: TOL=$10^{-6}$ and TOL=$10^{-3}$. \reviewerone{This distinction is motivated by the different applicable use cases. For example, a strict tolerance is often required when numerical simulations are used as a substitute for experiments, whereas a loose tolerance can be sufficient in scenarios such as solving inverse problems.}

The following test cases are shown: homogeneous waveguide (\S~\ref{sub:homogeneous_problem} and~\ref{sub:overlap}), open cavity (\S~\ref{sec:open_cavity_test}), Marmousi (\S~\ref{sub:marmousi}) and 3D Overthrust (\S~\ref{sub:overthrust}). \addedinrevision{The homogeneous waveguide and open cavity problems are discretized with P3 finite elements and 6 points per wavelength, the Marmousi test case is discretized with P2 elements and 8 points per wavelength, and the Overthrust problem is discretized with P1 elements and 10 points per wavelength. The implementation is done in the FreeFem++ domain specific language~\cite{Hecht:2012:NDF}. Our parallel implementation uses the MPI paradigm and we assign one core per subdomain. MUMPS~\cite{Amestoy:2001:FAM} solver is used to compute the local factorization in each subdomain.}

We also compare the use of zeroth-order ABC with that of Perfectly Matched Layers (PML) as interface conditions for the various sweeping algorithms. It is worth noting that in order to avoid transmitting attenuated solutions between subdomains, the width of the overlap region needs to be taken larger that the width of the PML. In the following numerical results, the overlap $\delta$ is equal to 4 mesh elements and the width of the PML is 2. The PML is defined following~\cite{bermudez2004exact}. \addedinrevision{In all tables, X means that the algorithm did not reach the convergence criterion in 2000 iterations for Jacobi, or in 400 iterations for the sweeping methods ; nppwl refers to the number of points per wavelength.} As we can see, using PML interface conditions, \addedinrevision{iteration counts are significantly better for all cases}.  Also note that as the number of subdomains increases, the superiority of PML becomes more pronounced.   
 
 When comparing the various algorithms, we see that the iteration counts are qualitatively in agreement with the spectral radius estimates summarized in Table~\ref{tab:algocvrate} but only when using PML as interface conditions. For the Jacobi method, the increase is linear w.r.t. to the number of subdomains. The sweeping methods have iteration counts that increase sublinearly with the number of subdomains. When using PML, the Block Jacobi (BJ) methods needs twice as many iterations as the other sweeping methods. We also see that \reviewertwo{the new preconditioner BSGS  offers no iteration count improvement over the BGS and SGS preconditioners although it is more expensive.}
  

\subsection{Homogeneous waveguide}
\label{sub:homogeneous_problem}

\begin{figure}[htp]
  \centering
  \includegraphics[width=0.8\textwidth]{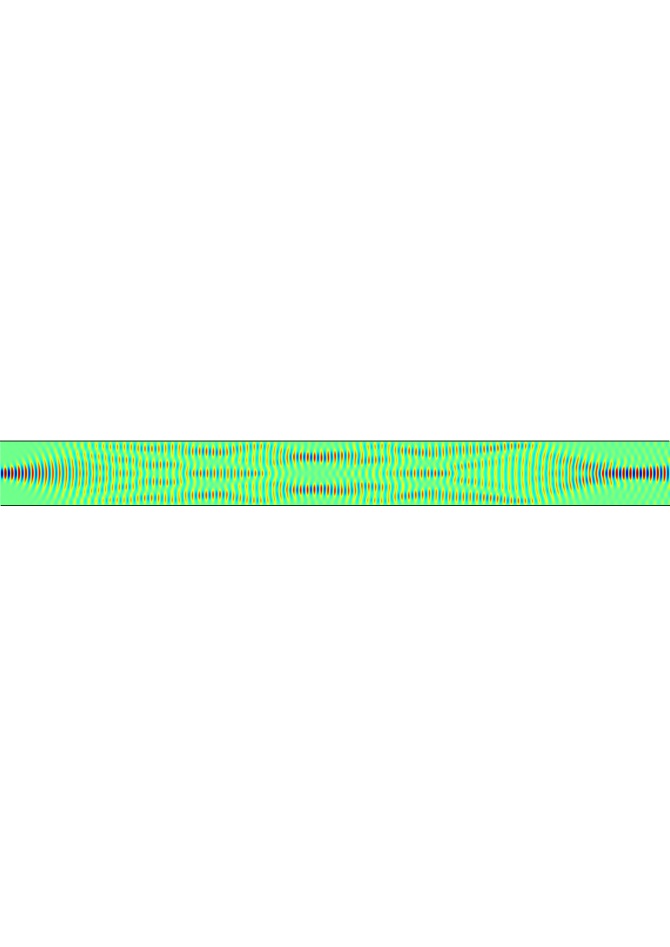}
  \caption{Homogeneous waveguide ($k=20\pi$)}
  \label{fig:HomogeneousSolution}
\end{figure}

First, we consider the homogeneous waveguide test case with a layered decomposition into $N$ subdomains. More specifically, we consider a rectangular geometry $\left( \Omega = [0, N] \times [0, 1] \right)$ made of a homogeneous medium. On the upper and lower sides of the waveguide, we impose homogeneous Dirichlet conditions (cf. black lines in Figure~\ref{fig:HomogeneousSolution}). In addition, we impose an absorbing boundary condition on the right side and \reviewerone{an inhomogeneous ABC on the left side where $g$ is the impedance operator applied to the excitation}. The global problem is written as
\begin{equation}
\label{homogeneous_problem}
\left\{
\begin{aligned}
\left(- k^2 - \Delta\right) u &= f\ \text{ in } \Omega\\
\left(\partial_{\vec{n}} + Ik \right) u &= 0 \ \text{ on }  \{x=N\}\times [0, 1] \\
\left(\partial_{\vec{n}} + Ik \right) u &= \reviewerone{g} \ \text{ on }  \{x=0\}\times [0, 1]\\
u &= 0 \ \text{ on } [0, N] \times \{y=0, y=1\}
\end{aligned}
\right.
\end{equation}
where $g =  e^{-120(y-0.5)^2}\sin(\pi y)$.

We considered two values for the wave number: $k =20$ and $k = 20\pi$ and results are given in Tables~\ref{tab:homowaveguide_vol_npplo6_P3_20} and~\ref{tab:homowaveguide_vol_npplo6_P3_20pi} for the zeroth order interface conditions and in Tables~\ref{tab:homowaveguide_vol_pmls_npplo6_P3_20} and~\ref{tab:homowaveguide_vol_pmls_npplo6_P3_20pi} for PML interface conditions. The iteration counts of the various sweeping algorithms are in agreement with the spectral radius estimates of Table~\ref{tab:algocvrate}. \addedinrevision{As $N$ increases, we can see that the iteration count mildly increases for PML, whereas it becomes much larger for the zeroth-order ABC, thereby demonstrating the superiority of PML over the zeroth-order ABC.}

\begin{table}[htp!]
\centering\begin{filecontents*}{fullsweepingwaveguide_ovl4_P3.csv}
f,p,PMLORAS,PMLBJS,PMLBGSS,PMLBSGSS,PMLSGSS,ROBINORAS,ROBINBJS,ROBINBGSS,ROBINBSGSS,ROBINSGSS,IPMLORAS,IPMLBJS,IPMLBGSS,IPMLBSGSS,IPMLSGSS,IROBINORAS,IROBINBJS,IROBINBGSS,IROBINBSGSS,IROBINSGSS
3.1831,5,27 ,14 ,13 ,14 ,14 ,29 ,10 ,5 ,4 ,5 ,14 ,6 ,7 ,8 ,8 ,13 ,4 ,2 ,2 ,2
10.0,5,29 ,9 ,5 ,4 ,4 ,33 ,13 ,9 ,8 ,9 ,9 ,2 ,1 ,1 ,1 ,14 ,4 ,3 ,3 ,2
3.1831,10,51 ,5 ,3 ,3 ,3 ,60 ,13 ,8 ,6 ,7 ,28 ,3 ,2 ,2 ,2 ,32 ,5 ,3 ,2 ,3
10.0,10,58 ,9 ,5 ,4 ,5 ,77 ,25 ,15 ,20 ,19 ,19 ,2 ,1 ,1 ,1 ,32 ,8 ,6 ,10 ,6
3.1831,20,109 ,6 ,3 ,3 ,3 ,122 ,20 ,11 ,13 ,11 ,58 ,4 ,2 ,2 ,2 ,68 ,8 ,5 ,4 ,5
10.0,20,122 ,9 ,6 ,4 ,5 ,173 ,44 ,30 ,46 ,38 ,39 ,2 ,1 ,1 ,1 ,67 ,15 ,11 ,20 ,9
3.1831,40,208 ,7 ,3 ,3 ,3 ,249 ,35 ,18 ,26 ,19 ,118 ,3 ,2 ,2 ,2 ,141 ,14 ,8 ,7 ,7
10.0,40,241 ,10 ,7 ,5 ,6 ,354 ,84 ,53 ,169 ,70 ,79 ,2 ,1 ,1 ,1 ,136 ,35 ,24 ,38 ,32
3.1831,80,417 ,6 ,3 ,3 ,3 ,507 ,56 ,32 ,45 ,35 ,238 ,3 ,2 ,2 ,2 ,282 ,23 ,10 ,14 ,10
10.0,80,462 ,11 ,7 ,6 ,6 ,703 ,153 ,104 ,X ,142 ,149 ,2 ,1 ,1 ,1 ,274 ,69 ,53 ,X ,69
\end{filecontents*}
\csvreader[tabular=|p{2cm}|p{2cm} p{2cm} p{2cm} p{2cm} p{2cm}|,table head=\hline N & Jacobi & BJ & BGS & BSGS & SGS\\\hline, table foot=\hline,head to column names,filter=\equal{\f}{3.1831}\and\not\equal{\p}{5}]{fullsweepingwaveguide_ovl4_P3.csv}{}
{\p & \ROBINORAS \ (\IROBINORAS) & \ROBINBJS \ (\IROBINBJS)& \ROBINBGSS \ (\IROBINBGSS) & \ROBINBSGSS \ (\IROBINBSGSS) &\ROBINSGSS \ (\IROBINSGSS)}
\caption{Homogeneous waveguide with zeroth-order ABC, $k=20$, $\delta = 4h$, TOL=$10^{-6} (10^{-3})$, nppwl = 6 , P3}
\label{tab:homowaveguide_vol_npplo6_P3_20}
\end{table}

\begin{table}[htp!]
\centering\csvreader[tabular=|p{2cm}|p{2cm} p{2cm} p{2cm} p{2cm} p{2cm}|,table head=\hline N & Jacobi & BJ & BGS & BSGS & SGS\\\hline, table foot=\hline,head to column names,filter=\equal{\f}{3.1831}\and\not\equal{\p}{5}]{fullsweepingwaveguide_ovl4_P3.csv}{}
{\p & \PMLORAS \ (\IPMLORAS) & \PMLBJS \ (\IPMLBJS)& \PMLBGSS \ (\IPMLBGSS) & \PMLBSGSS \ (\IPMLBSGSS) &\PMLSGSS \ (\IPMLSGSS)}
\caption{PML interface conditions, homogeneous waveguide, $k=20$, $\delta = 4h$, TOL=$10^{-6} (10^{-3})$, nppwl = 6 , P3 }
\label{tab:homowaveguide_vol_pmls_npplo6_P3_20}
\end{table}

\begin{table}[htp!]
\centering\csvreader[tabular=|p{2cm}|p{2cm} p{2cm} p{2cm} p{2cm} p{2cm}|,table head=\hline N & Jacobi & BJ & BGS & BSGS & SGS\\\hline, table foot=\hline,head to column names,filter=\equal{\f}{10.0}\and\not\equal{\p}{5}]{fullsweepingwaveguide_ovl4_P3.csv}{}
{\p & \ROBINORAS \ (\IROBINORAS) & \ROBINBJS \ (\IROBINBJS)& \ROBINBGSS \ (\IROBINBGSS) & \ROBINBSGSS \ (\IROBINBSGSS) &\ROBINSGSS \ (\IROBINSGSS)}
\caption{Homogeneous waveguide with zeroth-order ABC, $k=20\pi$, $\delta = 4h$, TOL=$10^{-6} (10^{-3})$, nppwl = 6 , P3}
\label{tab:homowaveguide_vol_npplo6_P3_20pi}
\end{table}

\begin{table}[htp!]
\centering\csvreader[tabular=|p{2cm}|p{2cm} p{2cm} p{2cm} p{2cm} p{2cm}|,table head=\hline N & Jacobi & BJ & BGS & BSGS & SGS\\\hline, table foot=\hline,head to column names,filter=\equal{\f}{10.0}\and\not\equal{\p}{5}]{fullsweepingwaveguide_ovl4_P3.csv}{}
{\p & \PMLORAS \ (\IPMLORAS) & \PMLBJS \ (\IPMLBJS)& \PMLBGSS \ (\IPMLBGSS) & \PMLBSGSS \ (\IPMLBSGSS) &\PMLSGSS \ (\IPMLSGSS)}
\caption{PML interface conditions, homogeneous waveguide, $k=20\pi$, $\delta = 4h$, TOL=$10^{-6} (10^{-3})$, nppwl = 6 , P3 }
\label{tab:homowaveguide_vol_pmls_npplo6_P3_20pi}
\end{table}

\subsection{Influence of the overlap}
\label{sub:overlap}

We have also tested the effect of the width of the overlap on the convergence in the case of the homogeneous waveguide. \addedinrevision{The overlap varies from 4 mesh sizes up to 16 mesh sizes. Interestingly, for zeroth-order ABC~(Table~\ref{tab:homowaveguide_vol_ovrl_N20_k20pi_P3}), we see that the number of iterations increases when going from $\delta=4$ to $\delta=8$, but decreases for $\delta=16$. For PML~(Table~\ref{tab:homowaveguide_vol_ovrl_N20_k20pi_P3_PML}), there is little if no effect of the overlap on the iteration count.}

\begin{table}[htp!]
\centering\begin{filecontents*}{fullsweepingwaveguideovl_P3.csv}
delt,f,p,PMLORAS,PMLBJS,PMLBGSS,PMLBSGSS,PMLSGSS,ROBINORAS,ROBINBJS,ROBINBGSS,ROBINBSGSS,ROBINSGSS,IPMLORAS,IPMLBJS,IPMLBGSS,IPMLBSGSS,IPMLSGSS,IROBINORAS,IROBINBJS,IROBINBGSS,IROBINBSGSS,IROBINSGSS
4, 3.1831, 20,109,6,3,3,3 , 122,20,11,13,11,58 ,4 ,2 ,2 ,2 ,68 ,8 ,5 ,4 ,5 
4, 10.0, 20,122,9,6,4,5 , 173,44,30,46,38,39 ,2 ,1 ,1 ,1 ,67 ,15 ,11 ,20 ,9
4, 3.1831, 40,208,7,3,3,3 , 249,35,18,26,19,118 ,3 ,2 ,2 ,2 ,141 ,14 ,8 ,7 ,7
4, 10.0, 40, 241,10,7,5,6 , 354,84,53,169,70 ,79 ,2 ,1 ,1 ,1 ,136 ,35 ,24 ,38 ,32
8, 3.1831, 20, 108,7,4,3,4 , 138,23,12,14,12,59 ,4 ,2 ,2 ,2 ,79 ,9 ,5 ,5 ,5 
8, 10.0, 20,123,8,5,4,4 , 168,54,27,69,36,39 ,2 ,1 ,1 ,1 ,70 ,19 ,6 ,7 ,6
8, 3.1831, 40,210,6,4,3,4 , 264,37,19,26,20,118 ,3 ,2 ,2 ,2 ,156 ,12 ,6 ,7 ,6
8, 10.0, 40,242,9,6,5,5 , 345,92,115,X,120,81 ,2 ,1 ,1 ,1 ,140 ,25 ,13 ,X ,16
16, 3.1831, 20,110,7,4,4,4 , 122,18,9,10,9,59 ,4 ,2 ,2 ,2 ,70 ,9 ,6 ,5 ,6 
16, 10.0, 20, 122,7,5,3,4 , 149,41,24,26,24,39 ,2 ,1 ,1 ,1 ,69 ,18 ,10 ,10 ,10
16, 3.1831, 40,214,6,4,3,4 , 235,28,15,16,15,118 ,3 ,2 ,2 ,2 ,130 ,11 ,6 ,6 ,6
16, 10.0, 40, 241,8,5,4,5 , 306,72,44,49,45,81 ,2 ,1 ,1 ,1 ,145 ,33 ,20 ,25 ,20
\end{filecontents*}
\csvreader[tabular=|p{2cm}|p{2cm} p{2cm} p{2cm} p{2cm} p{2cm}|,table head=\hline $\delta$ & Jacobi & BJ & BGS & BSGS & SGS\\\hline, table foot=\hline,head to column names,filter=\equal{\f}{10.0}\and\equal{\p}{20}]{fullsweepingwaveguideovl_P3.csv}{}
{\delt & \ROBINORAS \ (\IROBINORAS) & \ROBINBJS \ (\IROBINBJS)& \ROBINBGSS \ (\IROBINBGSS) & \ROBINBSGSS \ (\IROBINBSGSS) &\ROBINSGSS \ (\IROBINSGSS)}
\caption{Homogeneous waveguide with zeroth-order ABC, $k=20\pi$, $N = 20$, $\delta$ varies, TOL=$10^{-6} (10^{-3})$, nppwl = 6 , P3 }
\label{tab:homowaveguide_vol_ovrl_N20_k20pi_P3}
\end{table}

\begin{table}[htp!]
\centering\csvreader[tabular=|p{2cm}|p{2cm} p{2cm} p{2cm} p{2cm} p{2cm}|,table head=\hline $\delta$ & Jacobi & BJ & BGS & BSGS & SGS\\\hline, table foot=\hline,head to column names,filter=\equal{\f}{10.0}\and\equal{\p}{20}]{fullsweepingwaveguideovl_P3.csv}{}
{\delt & \PMLORAS \ (\IPMLORAS) & \PMLBJS \ (\IPMLBJS)& \PMLBGSS \ (\IPMLBGSS) & \PMLBSGSS \ (\IPMLBSGSS) &\PMLSGSS \ (\IPMLSGSS)}
\caption{Homogeneous waveguide with PML interface conditions, $k=20\pi$, $N = 20$, $\delta$ varies, TOL=$10^{-6} (10^{-3})$, nppwl = 6 , P3 }
\label{tab:homowaveguide_vol_ovrl_N20_k20pi_P3_PML}
\end{table}

%

\subsection{Open cavity test}
  \label{sec:open_cavity_test}

\begin{figure}[htp!]
  \centering
  \includegraphics[width=0.8\textwidth]{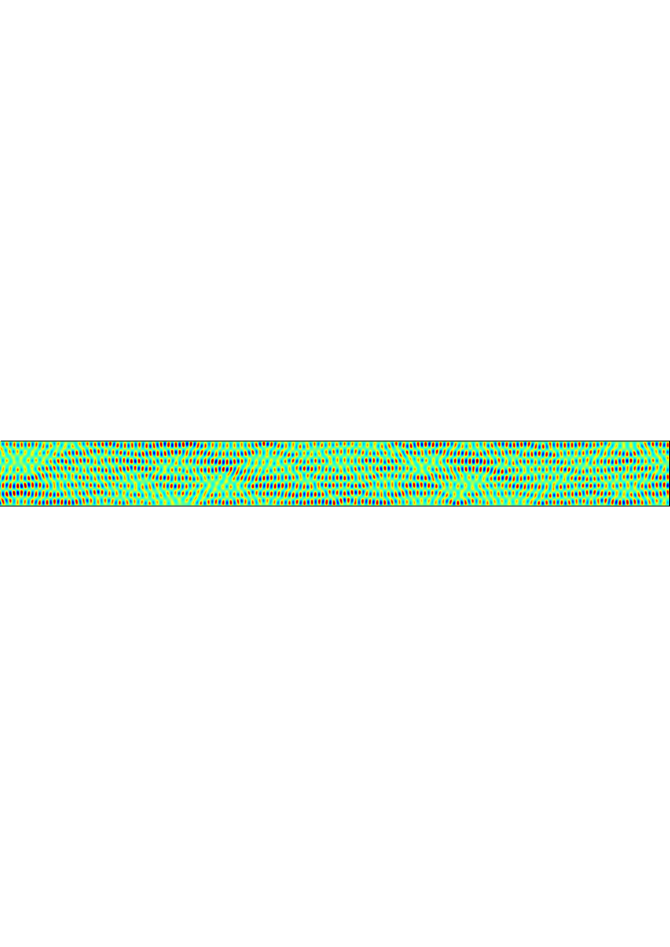}
  \caption{Open cavity solution ($k=20\pi$)}
  \label{fig:CavitySolution}
\end{figure}

Same as before, the domain is rectangular with an homogeneous medium and its length increases with the number of subdomains. The open cavity test is challenging due to the homogeneous Dirichlet conditions imposed on three sides (cf. black lines in Figure~\ref{fig:CavitySolution}). In addition, we perform an excitation on the left side \reviewerone{which is modelled by the right hand side $g$ which is the impedance operator applied to the excitation}. The Dirichlet conditions create rebounds leading to an increase in the number of reflections, this phenomenon is exacerbated for high-frequency regimes. The global problem can be written as
\begin{equation}
\left\{
\begin{aligned}
\left(- k^2 - \Delta\right) u &= f\ \text{ in } \Omega\\
\left(\partial_{\vec{n}} + Ik \right) u &= g \ \text{ on }  \Gamma \\
u &= 0 \ \text{ on } \partial\Omega \backslash \Gamma
\end{aligned}
\right.
\end{equation}
where $\Gamma\vcentcolon = \{x=0\}\times [0,1]$ and $g =  \exp^{-ik\left(x \cos(\theta) + y \sin(\theta)\right)},\ \theta = \frac{\pi}{8}$. It corresponds to an incident plane wave propagating at an angle $\frac{\pi}{8}$ w.r.t. to the horizontal direction, see Fig.~\ref{fig:CavitySolution}. This creates numerous reflections on the lateral boundaries of the open cavity. Here as well, there is roughly a factor two in the iteration counts in favour of the SGS, BGS and BSGS algorithms compared to BJ algorithm.

\addedinrevision{Once again, we see that PML outperforms zeroth-order ABC, although we observe a significantly higher increase in the iteration count when transitioning from $k = 20$ to $k=20\pi$, compared to the homogeneous waveguide. This is most likely due to the fact that a wavenumber $k=20\pi$ creates more trapping and resonance effects.}

\begin{table}[htp!]
\centering\begin{filecontents*}{fullsweepingcavity_ovl4_P3.csv}
f,p,PMLORAS,PMLBJS,PMLBGSS,PMLBSGSS,PMLSGSS,ROBINORAS,ROBINBJS,ROBINBGSS,ROBINBSGSS,ROBINSGSS,IPMLORAS,IPMLBJS,IPMLBGSS,IPMLBSGSS,IPMLSGSS,IROBINORAS,IROBINBJS,IROBINBGSS,IROBINBSGSS,IROBINSGSS
3.1831,5,43 ,19 ,16 ,18 ,17 ,47 ,14 ,8 ,8 ,9 ,26 ,10 ,9 ,10 ,10 ,31 ,10 ,5 ,5 ,5
10.0,5,75 ,24 ,14 ,12 ,12 ,80 ,32 ,20 ,17 ,19 ,29 ,8 ,4 ,4 ,4 ,34 ,12 ,8 ,9 ,9
3.1831,10,93 ,12 ,6 ,6 ,6 ,103 ,20 ,10 ,11 ,11 ,56 ,8 ,4 ,4 ,4 ,64 ,13 ,7 ,7 ,7
10.0,10,169 ,30 ,19 ,16 ,17 ,191 ,57 ,35 ,34 ,37 ,56 ,8 ,4 ,4 ,4 ,81 ,23 ,15 ,20 ,19
3.1831,20,193 ,12 ,6 ,6 ,6 ,217 ,30 ,15 ,18 ,16 ,117 ,8 ,4 ,4 ,4 ,139 ,21 ,11 ,14 ,11
10.0,20,329 ,39 ,28 ,23 ,24 ,388 ,98 ,62 ,70 ,68 ,117 ,9 ,5 ,5 ,5 ,193 ,45 ,29 ,62 ,37
3.1831,40,398 ,12 ,6 ,6 ,6 ,440 ,49 ,24 ,31 ,26 ,241 ,8 ,4 ,4 ,4 ,293 ,36 ,18 ,27 ,20
10.0,40,562 ,40 ,31 ,25 ,27 ,690 ,142 ,90 ,175 ,110 ,209 ,8 ,5 ,4 ,4 ,365 ,80 ,54 ,170 ,72
3.1831,80,797 ,12 ,6 ,6 ,6 ,875 ,78 ,39 ,66 ,44 ,375 ,7 ,3 ,3 ,3 ,533 ,60 ,30 ,62 ,32
10.0,80,1115 ,41 ,32 ,26 ,27 ,1330 ,226 ,245 ,X ,251 ,356 ,8 ,5 ,4 ,5 ,728 ,154 ,106 ,X ,146
\end{filecontents*}
\csvreader[tabular=|p{2cm}|p{2cm} p{2cm} p{2cm} p{2cm} p{2cm}|,table head=\hline N & Jacobi & BJ & BGS & BSGS & SGS\\\hline, table foot=\hline,head to column names,filter=\equal{\f}{3.1831}\and\not\equal{\p}{5}]{fullsweepingcavity_ovl4_P3.csv}{}
{\p & \ROBINORAS \ (\IROBINORAS) & \ROBINBJS \ (\IROBINBJS)& \ROBINBGSS \ (\IROBINBGSS) & \ROBINBSGSS \ (\IROBINBSGSS) &\ROBINSGSS \ (\IROBINSGSS)}
\caption{Zeroth-order ABC, open cavity, $k=20$, $\delta = 4h$, TOL=$10^{-6} (10^{-3})$, nppwl = 6 , P3 }
\label{tab:opencavity_vol_npplo6_P3_20}
\end{table}

\begin{table}[htp!]
\centering\csvreader[tabular=|p{2cm}|p{2cm} p{2cm} p{2cm} p{2cm} p{2cm}|,table head=\hline N & Jacobi & BJ & BGS & BSGS & SGS\\\hline, table foot=\hline,head to column names,filter=\equal{\f}{3.1831}\and\not\equal{\p}{5}]{fullsweepingcavity_ovl4_P3.csv}{}
{\p & \PMLORAS \ (\IPMLORAS) & \PMLBJS \ (\IPMLBJS)& \PMLBGSS \ (\IPMLBGSS) & \PMLBSGSS \ (\IPMLBSGSS) &\PMLSGSS \ (\IPMLSGSS)}
\caption{PML interface conditions, open cavity, $k=20$, $\delta = 4h$, TOL=$10^{-6} (10^{-3})$, nppwl = 6 , P3 }
\label{tab:opencavity_vol_pmls_npplo6_P3_20}
\end{table}

\begin{table}[htp!]
\centering\csvreader[tabular=|p{2cm}|p{2cm} p{2cm} p{2cm} p{2cm} p{2cm}|,table head=\hline N & Jacobi & BJ & BGS & BSGS & SGS\\\hline, table foot=\hline,head to column names,filter=\equal{\f}{10.0}\and\not\equal{\p}{5}]{fullsweepingcavity_ovl4_P3.csv}{}
{\p & \ROBINORAS \ (\IROBINORAS) & \ROBINBJS \ (\IROBINBJS)& \ROBINBGSS \ (\IROBINBGSS) & \ROBINBSGSS \ (\IROBINBSGSS) &\ROBINSGSS \ (\IROBINSGSS)}
\caption{Zeroth-order ABC, open cavity, $k=20\pi$, $\delta = 4h$, TOL=$10^{-6} (10^{-3})$, nppwl = 6 , P3 }
\label{tab:opencavity_vol_npplo6_P3_20pi}
\end{table}

\begin{table}[htp!]
\centering\csvreader[tabular=|p{2cm}|p{2cm} p{2cm} p{2cm} p{2cm} p{2cm}|,table head=\hline N & Jacobi & BJ & BGS & BSGS & SGS\\\hline, table foot=\hline,head to column names,filter=\equal{\f}{10.0}\and\not\equal{\p}{5}]{fullsweepingcavity_ovl4_P3.csv}{}
{\p & \PMLORAS \ (\IPMLORAS) & \PMLBJS \ (\IPMLBJS)& \PMLBGSS \ (\IPMLBGSS) & \PMLBSGSS \ (\IPMLBSGSS) &\PMLSGSS \ (\IPMLSGSS)}
\caption{PML interface conditions, open cavity, $k=20\pi$, $\delta = 4h$, TOL=$10^{-6} (10^{-3})$, nppwl = 6 , P3 }
\label{tab:opencavity_vol_pmls_npplo6_P3_20pi}
\end{table}

\subsection{Marmousi} 
\label{sub:marmousi}

\begin{figure}[htp!]
\centering
   \includegraphics[width=0.9\textwidth, clip=true, trim = 0cm 15cm 0cm 15cm]{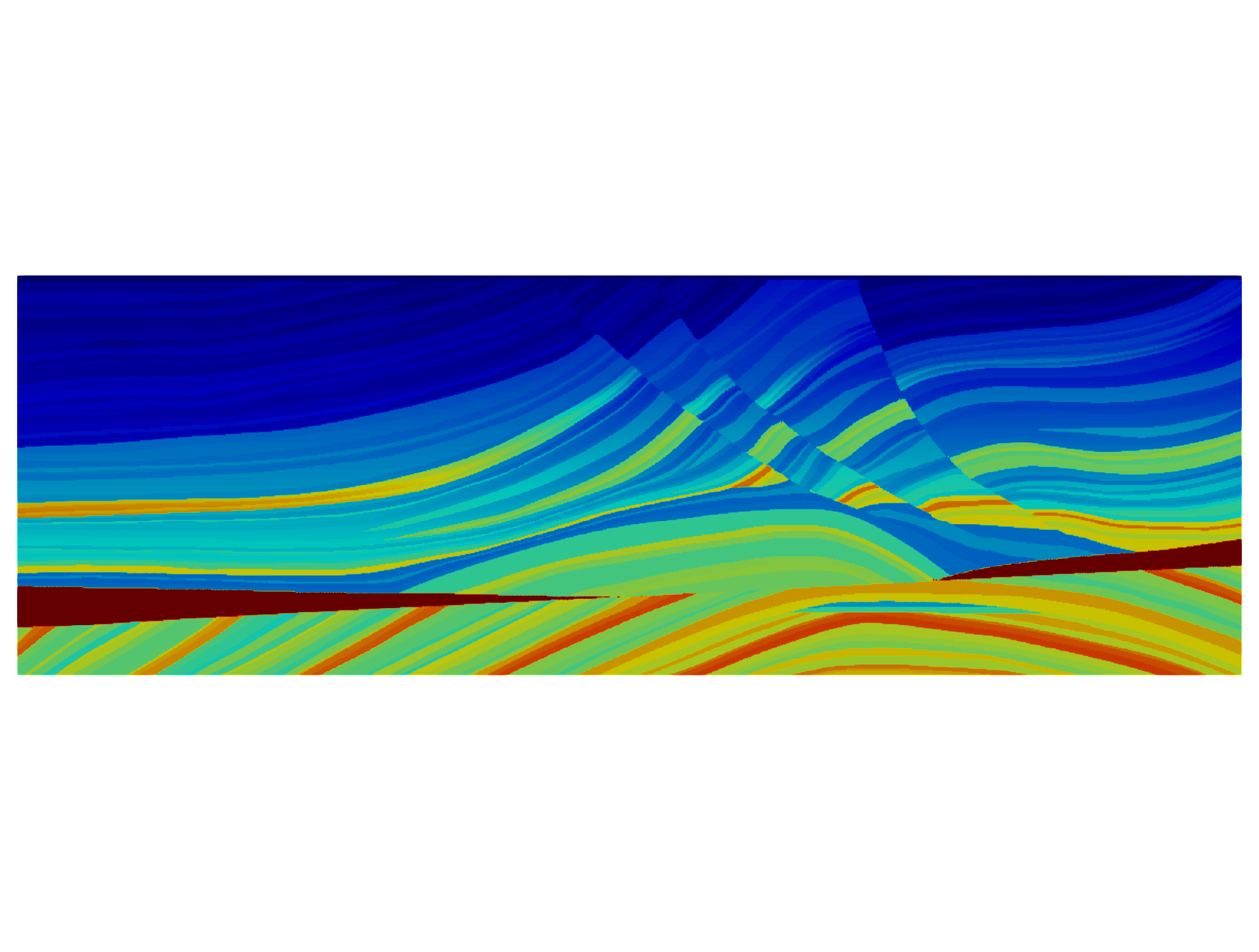}
   \caption{Velocity model of the Marmousi test case}
   \label{fig:marmousi}
\end{figure}

\begin{figure}[htp!]
\centering
	\includegraphics[width=0.9\textwidth, clip=true, trim = 0cm 30cm 0cm 30cm]{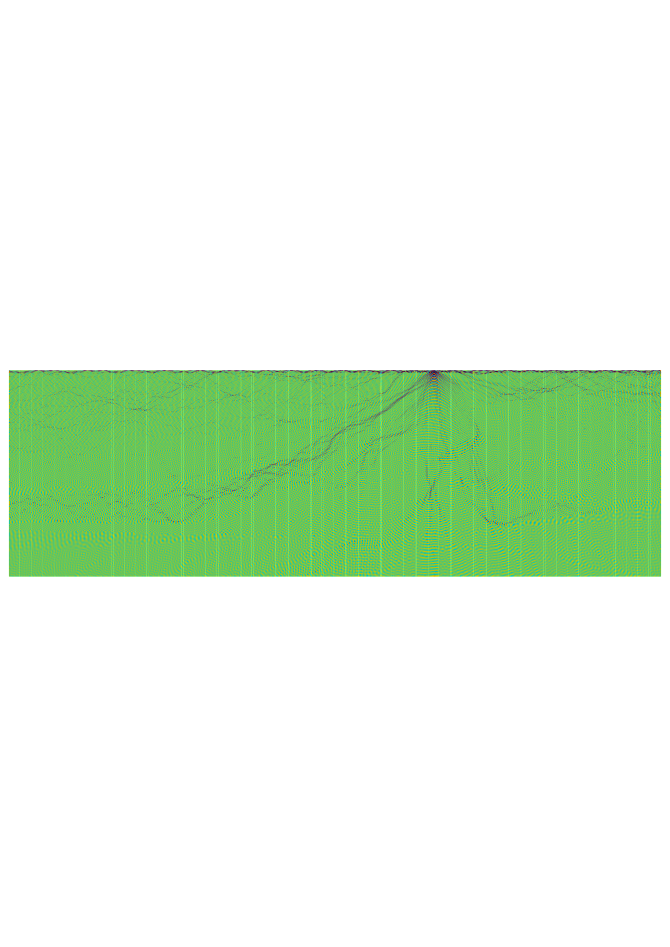}
	\caption{Real part of the solution for $f = 100$ Hz for the Marmousi test case}
	\label{fig:marmousi_f100}
\end{figure}

Here we consider the Marmousi benchmark~\cite{Brougois:1990:marmousi}. The velocity profile of the model is depicted in Figure~\ref{fig:marmousi}. The domain is of size 9.2 km $\times$ 3 km. A Neumann boundary condition is imposed at the top boundary, and PMLs are used on the other three boundaries. The source is located at the top. The problem is discretized with P2 elements on a regular mesh with 8 points per wavelength for the reference wavelength corresponding to $c = 2$. Figure~\ref{fig:marmousi_f100} shows the real part of the acoustic field at 100 Hz frequency. We present results for 25 Hz, 50 Hz and 100 Hz frequencies using the volumic preconditioner. The number of degrees of freedom grows from 1.1 million for $f = 25$ Hz to 4.5 million for $f=50$ Hz and 17.8 million for $f = 100$ Hz. Table~\ref{tab:marm_vol_robin_npplo8_P2} shows numerical results for zeroth order ABC interface conditions and Table~\ref{tab:marm_vol_pmls_npplo8_P2} for PML. Interestingly we see that for ABC results are quite bad with high iteration counts and dependence on both the number of subdomains and the wave number. Whereas with PML interface conditions, both the iteration counts and their behaviour with respect to the number of subdomains and wavenumber is similar to the previous test cases on waveguides and open cavity. This very bad behaviour of zeroth order ABC at higher frequency and subdomain count can be explained by the variability of the coefficients. 

\begin{table}[htp!]
\begin{center}

\begin{tabular}{ |p{2cm} | p{2cm}  p{2cm}  p{2cm} p{2cm} p{2cm} |}
 \hline
  & \multicolumn{5}{c|}{25 Hz}\\
 \hline
 N & Jacobi & BJ & BGS & BSGS & SGS \\
\hline
 3   & 22 (12)& 11 (6)& 10 (5) & 7 (4) & 7 (4) \\
 7   & 46 (25)& 17 (9)& 13 (7) & 10 (6) & 11 (6) \\
 14  & 94 (51)&25 (14)&18 (10)&19 (12)&20(12) \\
 28  & 185 (100)&41 (24)&27 (16)&30 (20)&30 (19) \\
 56  & 382 (220)&98 (69)&65 (47)&101 (91)&70 (49) \\
 \hline
\end{tabular}
\begin{tabular}{ |p{2cm} | p{2cm}  p{2cm}  p{2cm} p{2cm} p{2cm} |}
 \hline
  & \multicolumn{5}{c|}{50 Hz}\\
 \hline
 N & Jacobi & BJ & BGS & BSGS & SGS \\
\hline
 7   & 49 (25) & 19 (9) & 16 (8) & 13 (8) & 14 (8) \\
 14  & 98 (51)&28 (15)&22 (12)&21 (13)&21 (11) \\
 28  & 195 (101)&49 (26)&38 (21)&60 (48)&47 (32)  \\
 56  & 426 (222)&123(81)&90 (61)& X (192)&143 (116) \\
 112 & 1505 (690)& X (324)& X (X)& X (X)& X (X) \\
 \hline
\end{tabular}
\begin{tabular}{ |p{2cm} | p{2cm}  p{2cm}  p{2cm} p{2cm} p{2cm} |}
 \hline
  & \multicolumn{5}{c|}{100 Hz}\\
 \hline
 N & Jacobi & BJ & BGS & BSGS & SGS \\
\hline
 56  & 476 (216)& 147 (87)& 120 (77) & X (X) & 208 (171) \\
 112 & 1701 (691)& X (X)& X (375) & X (X) & X (X) \\
 224 & X (X)& X (X)& X (X) & X (X) & X (X) \\
 \hline
\end{tabular}

\end{center}
\caption{Zeroth-order ABC, Marmousi test case, $\delta = 4h$, TOL=$10^{-6} (10^{-3})$, nppwl = 8, P2. X means that the algorithm did not reach the convergence criterion in 2000 iterations for Jacobi, or in 400 iterations for the sweeping methods.}
\label{tab:marm_vol_robin_npplo8_P2}
\end{table}


\begin{table}[htp!]

\begin{center}
\begin{tabular}{ |p{2cm} | p{2cm}  p{2cm}  p{2cm} p{2cm} p{2cm} |}
 \hline
  & \multicolumn{5}{c|}{25 Hz}\\
 \hline
 N & Jacobi & BJ & BGS & BSGS & SGS \\
\hline
 3   & 14 (8)& 7 (4)& 6 (3) & 4 (2) & 5 (3) \\
 7   & 33 (17)& 10 (5)&8 (4)&5 (3)&6 (4) \\
 14  & 64 (35)& 11 (6)&8 (5)&6 (4)&7 (4) \\
 28  & 126 (66)&13 (7)&9 (5)&7 (4)&8 (4) \\
 56  & 247 (124)&18 (11)&12 (7)&11 (6)&12 (6) \\
 \hline
\end{tabular}
\begin{tabular}{ |p{2cm} | p{2cm}  p{2cm}  p{2cm} p{2cm} p{2cm} |}
 \hline
  & \multicolumn{5}{c|}{50 Hz}\\
 \hline
 N & Jacobi & BJ & BGS & BSGS & SGS \\
\hline
 7   & 34 (18) &10 (6)&9 (4)&6 (3)&6 (4)  \\
 14  & 69 (35)&12 (6)&10 (5)&7 (4)&8 (4) \\
 28  & 133 (66)&14 (7)&11 (6)&7 (4)&8 (5) \\
 56  & 260 (125)&18 (9)&13 (7)&10 (5)&11 (5) \\
 112 & 531 (240)&32 (17)&22 (12)&20 (11)& 21 (11)\\
 \hline
\end{tabular}
\begin{tabular}{ |p{2cm} | p{2cm}  p{2cm}  p{2cm} p{2cm} p{2cm} |}
 \hline
  & \multicolumn{5}{c|}{100 Hz}\\
 \hline
 N & Jacobi & BJ & BGS & BSGS & SGS \\
\hline
 56  & 242 (103)& 18 (9)& 14 (6) & 10 (5) & 11 (5) \\
 112 & 483 (198)& 29 (13)& 20 (9) & 15 (7) & 18 (8) \\
 224 & 1069 (417)& 66 (32)& 42 (20) & 62 (44) & 38 (18) \\
 \hline
\end{tabular}
\end{center}

\caption{PML interface conditions, Marmousi test case, $\delta = 4h$, TOL=$10^{-6} (10^{-3})$, nppwl = 8, P2 }
\label{tab:marm_vol_pmls_npplo8_P2}
\end{table}


\subsection{3D Overthrust} 
\label{sub:overthrust}

\begin{figure}[htp!]
\centering
	\includegraphics[width=0.7\textwidth, clip=true, trim = 0cm 0cm 0cm 0cm]{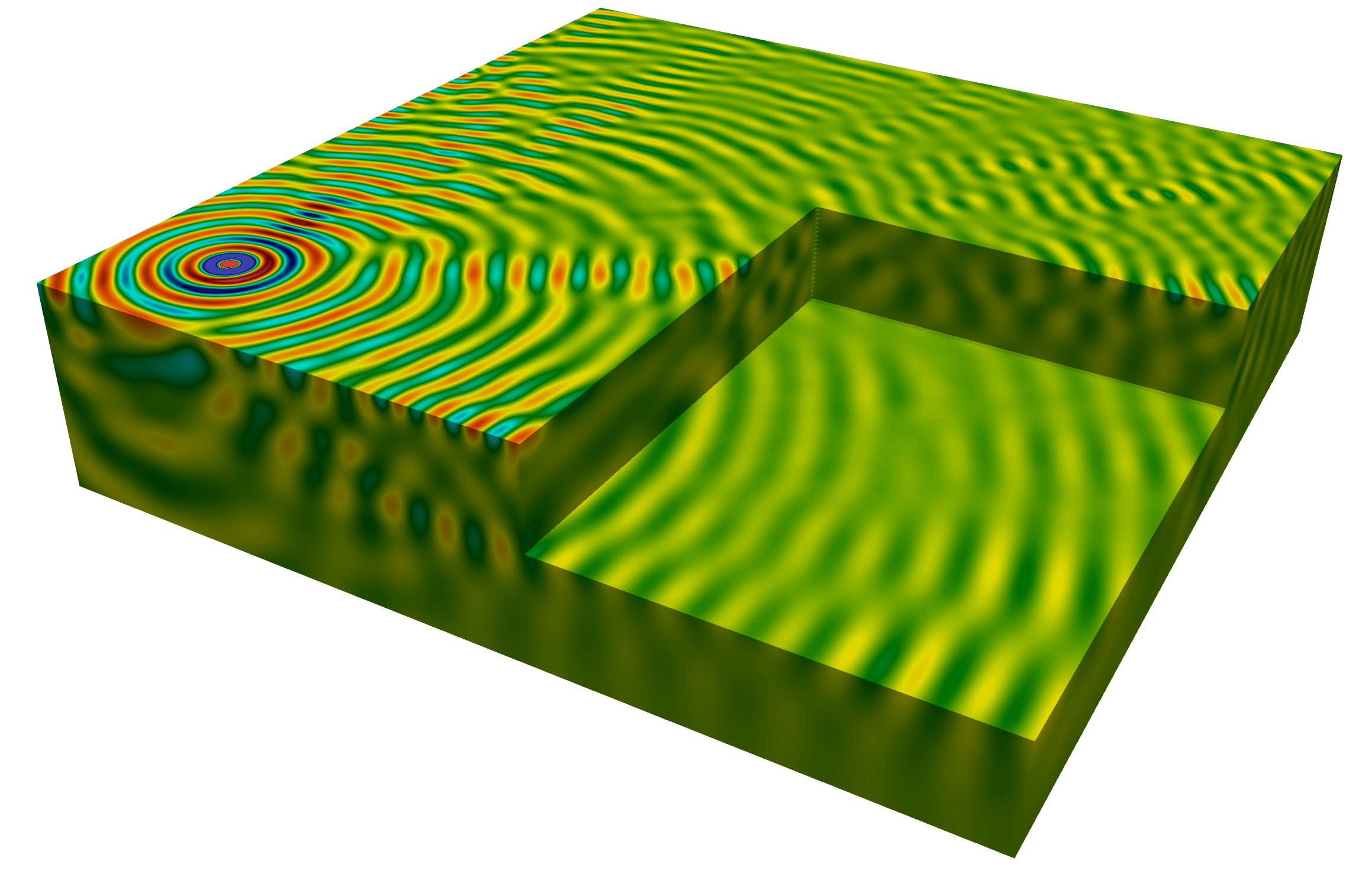}
	\caption{Real part of the solution for $f = 4$ Hz for the 3D Overthrust benchmark}
	\label{fig:overthrust_f4}
\end{figure}

Here we consider the 3D Overthrust acoustic benchmark. The velocity profile of the model is depicted in Figure~\ref{fig:marmousi}. The domain is of size 20 km $\times$ 20 km $\times$ 4.65 km. A Neumann boundary condition is imposed at the top boundary, and PMLs are used on the other five boundaries. The source is located at the top, at (2.5 km, 2.5 km). Figure~\ref{fig:overthrust_f4} shows the real part of the acoustic field at 4 Hz frequency. The problem is discretized on a regular mesh with 10 points per wavelength for the reference wavelength corresponding to $c = 2$. \reviewerone{Due to memory constraints, we were not able to use higher order discretizations than P1 elements. Indeed, sweeping methods are more constrained than standard domain decomposition methods in terms of subdomains size since the minimal size of a subdomain corresponds to $n_{d}\delta$, where $n_{d}$ is the number of d.o.f.s in a slice, and $\delta$ the width of the overlap. A possible fix could be to use more than one core per subdomain, in particular by considering a recursive method using again a sweeping preconditioner in the subdomains.}

Tables~\ref{tab:ovt_vol_robin_npplo10_P1} and~\ref{tab:ovt_vol_pmls_npplo10_P1} present results for 1 Hz, 2 Hz and 4 Hz frequencies using the volumic preconditioner with zeroth order ABC and PML interface conditions. The number of degrees of freedom grows from 0.33 million for $f = 1$ Hz to 16.2 million for $f = 4$ Hz. We observe similar trends as for the Marmousi test case, although the zeroth order ABC case is not so bad, as the frequency and number of wavelengths in one direction is not as high for this 3D test case. Moreover, we can see that the iteration count is much more dependent on the number of subdomains than on the frequency.

\begin{figure}[htp!]
\centering
   \includegraphics[width=0.7\textwidth, clip=true, trim = 0cm 1.2cm 0cm 1.2cm]{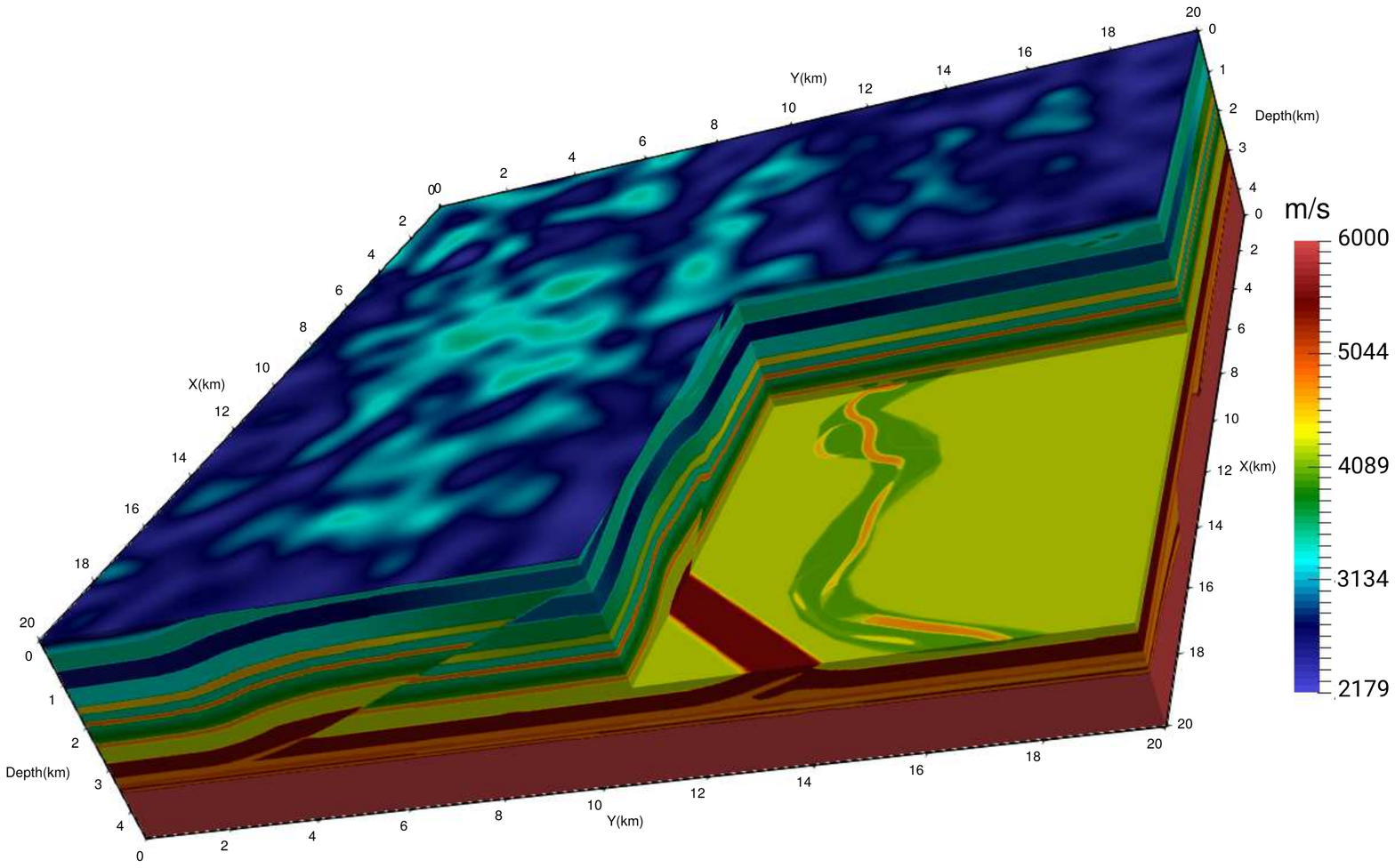}
   \caption{Velocity model of the 3D Overthrust benchmark}
   \label{fig:overthrust}
\end{figure}

\begin{table}[htp!]
\begin{center}

\begin{tabular}{ |p{2cm} | p{2cm}  p{2cm}  p{2cm} p{2cm} p{2cm} |}
 \hline
  & \multicolumn{5}{c|}{1 Hz} \\
 \hline
 N & Jacobi & BJ & BGS & BSGS & SGS \\
\hline
 3   & 14 (7)& 7 (4)& 6 (3) & 4 (2) & 4 (2) \\
 7   & 32 (17)& 9 (5)& 7 (4)& 5 (3)&6 (3) \\
 14  & 94 (45)& 18 (9)&10 (5)&9 (5)&9 (5) \\
 \hline
\end{tabular}
\begin{tabular}{ |p{2cm} | p{2cm}  p{2cm}  p{2cm} p{2cm} p{2cm} |}
 \hline
  & \multicolumn{5}{c|}{2 Hz} \\
 \hline
 N & Jacobi & BJ & BGS & BSGS & SGS \\
\hline
 7   & 34(17) &10 (5)&8 (4)&5 (3)&6 (3)  \\
 14  & 68 (34)&15 (8)&10 (5)&8 (4)&9 (4) \\
 28  & 163 (80)&30 (16)&16 (9)&17 (10)&14 (8) \\
 \hline
\end{tabular}	
\begin{tabular}{ |p{2cm} | p{2cm}  p{2cm}  p{2cm} p{2cm} p{2cm} |}
 \hline
  & \multicolumn{5}{c|}{4 Hz} \\
 \hline
 N & Jacobi & BJ & BGS & BSGS & SGS \\
\hline
28  &162 (80)&32 (17)&19 (10)&24 (15)&17 (9) \\
56  &365 (187)&87 (48)&42 (22)&168 (133)&36 (20) \\
 \hline
\end{tabular}
\end{center}
\caption{Zeroth-order ABC, Overthrust test case, $\delta = 4h$, TOL=$10^{-6} (10^{-3})$, nppwl = 10 , P1 }
\label{tab:ovt_vol_robin_npplo10_P1}
\end{table}

\begin{table}[htp!]
\begin{center}

\begin{tabular}{ |p{2cm} | p{2cm}  p{2cm}  p{2cm} p{2cm} p{2cm} |}
 \hline
  & \multicolumn{5}{c|}{1 Hz} \\
 \hline
 N & Jacobi & BJ & BGS & BSGS & SGS \\
\hline
 3   & 11 (7)& 5 (3)& 4 (2) & 3 (2) & 3 (2) \\
 7   & 28 (15)& 7 (4)&5 (3)&4 (2)&5 (3) \\
 14  & 59 (31)& 9 (5)&6 (4)&5 (3)&6 (3) \\
 \hline
\end{tabular}
\begin{tabular}{ |p{2cm} | p{2cm}  p{2cm}  p{2cm} p{2cm} p{2cm} |}
 \hline
  & \multicolumn{5}{c|}{2 Hz} \\
 \hline
 N & Jacobi & BJ & BGS & BSGS & SGS \\
\hline
 7   & 30(16) &8 (4)&6 (3)&4 (2)&5 (3)  \\
 14  & 63 (31)&11 (6)&7 (4)&6 (3)&7 (3) \\
 28  & 123 (61)&16 (8)&10 (5)&7 (4)&9 (5) \\
 \hline
\end{tabular}
\begin{tabular}{ |p{2cm} | p{2cm}  p{2cm}  p{2cm} p{2cm} p{2cm} |}
 \hline
  & \multicolumn{5}{c|}{4 Hz} \\
 \hline
 N & Jacobi & BJ & BGS & BSGS & SGS \\
\hline
28  &121 (63)&18 (9)&11 (6)&10 (6)&11 (6) \\
56  &224 (122)&28 (15)&18 (9)&19 (12)&16 (9) \\
 \hline
\end{tabular}
\end{center}

\caption{PML interface conditions, Overthrust test case, $\delta = 4h$, TOL=$10^{-6} (10^{-3})$, nppwl = 10 , P1 }
\label{tab:ovt_vol_pmls_npplo10_P1}
\end{table}

\section{Conclusion} 
\label{sec:conclusion}

We have introduced a unified framework for several sweeping algorithms which eases their derivation and theoretical studies, see Table~\ref{tab:algocvrate}. We have focused on three double sweep methods: SGS of~\cite{Nataf:1997:CRO,vion:2018:improved}, \reviewertwo{BGS of~\cite{Stolk:2013:rapidly} and BJ of~\cite{vion:2014:double}}. We conclude both from their theoretical convergence rates and from our numerical results that the SGS \reviewertwo{and BGS} algorithms have a better behavior than the BJ algorithm. \reviewertwo{Note also that the new algorithm BSGS is not competitive with respect to the other three preconditioners.} The numerical experiments also illustrate the superiority of PML as interface conditions compared to zeroth order ABCs, especially in the variable coefficient case. Although these algorithms are more easily derived in their substructured form, they lend themselves to volumic variants that allow for inexact solves in the subdomains, see \S~\ref{sub:volumicalgo}. In this form, these methods should have a behavior similar to the sweeping algorithm of~\cite{Engquist:2011:sweeping}. \newline

Also, an intrinsic problem with double sweep methods is that, due to the sequentiality of the algorithm, subdomains are idle most of the time. To be more precise, in SGS and BGS preconditioners, only one process is active at a time, whereas in BJ, two processes are active simultaneously (one for the left sweep and one for the right sweep). One way to overcome this limitation, in the context of multiple right-hand side problems, is to introduce a pipelining technique to improve parallelism, resulting in significant speed-ups, see~\cite{Stolk:2013:rapidly,vion:2014:double}.\newline

We have considered stripwise decompositions and it would be interesting to see if the present analysis can be extended to L-sweeps preconditioners (see e.g.~\cite{Taus:2020:Lsweeps,dai:2022:multidirectional} and references therein) used for checkerboard type decompositions.

\section*{Acknowledgment}

This work was granted access to the HPC resources of OCCIGEN@CINES under the allocation  2021-067730 granted by GENCI. The first author was partially supported by TOTAL Centre  Scientifique  et  Technique  Jean  F\'eger.

\bibliographystyle{plain}
\bibliography{bookddm}

\end{document}